\magnification=\magstep1
\input amstex
\voffset=-3pc
\documentstyle{amsppt}

\vsize=9.5truein
\hsize=6.5truein

\loadeufm
\loadmsbm
\def\sH{\Cal H}
\def\sS{\Cal S}

\def\sU{\Cal U}
\def\sI{\Cal I}
\def\cV{\Cal V}
\def\sB{\Cal B}
\def\bC{\Bbb C}
\def\bR{\Bbb R}
\def\sK{\Cal K}
\magnification=\magstep1
\parskip=6pt

\NoBlackBoxes
\topmatter
\title  Close Hereditary $C^*$-Subalgebras and the Structure of Quasi-Multipliers\endtitle
\rightheadtext{Close Hereditary $C^*$-subalgebras}
\leftheadtext{Lawrence G.~Brown}

\author Lawrence G.~Brown$^*$
\endauthor
\abstract{We answer a question of Takesaki by showing that the following can be derived from the thesis of N-T Shen:  If $A$ and $B$ are $\sigma$-unital hereditary $C^*$-subalgebras of $C$ such that $||p-q||<1$, where $p$ and $q$ are the corresponding open projections, then $A$ and $B$ are isomorphic.  We give some further elaborations and counterexamples with regard to the $\sigma$-unitality hypothesis.  We produce a natural one-to-one correspondence between complete order isomorphisms of $C^*$-algebras and invertible left multipliers of imprimitivity bimodules.  A corollary of the above two results is that any complete order isomorphism between $\sigma$-unital $C^*$-algebras is the composite of an isomorphism with an inner complete order isomorphism.  We give a separable counterexample to a question of Akemann and Pedersen; namely, the space of quasi-multipliers is not linearly generated by left and right multipliers.  But we show that the space of quasi-multipliers is multiplicatively generated by left and right multipliers in the $\sigma$-unital case.  In particular every positive quasi-multiplier is of the form $T^*T$ for $T$ a left multiplier.  We show that a Lie theory consequence of the negative result just stated is that the map sending $T$ to $T^*T$ need not be open, even for very nice $C^*$-algebras.  We show that surjective maps between $\sigma$-unital $C^*$-algebras induce surjective maps on left, right, and quasi-multipliers.  (The more significant similar result for multipliers is Pedersen's non-commutative Tietze extension theorem.)  We elaborate the relations of the above with continuous fields of Hilbert spaces and in so doing answer a question of Dixmier and Douady (yes for separable fields, no in general).  We discuss the relationship of our results to the theory of perturbations of $C^*$-algebras.}
\endabstract
\endtopmatter
$*$\ This work was done while the author was visiting the Mathematical Sciences Research Institute and was partially supported by M.S.R.I and the National Science Foundation.

\subhead \S1.\ Introduction\endsubhead

In her thesis, N.T.~Shen gave a characterization of the relative position of two hereditary $C^*$--subalgebras of a $C^*$--algebra.  (An alternative proof of Shen's main result was given in [9].)
Recently, M.~Takesaki asked whether Shen's result could be used to answer the following:\ If $p$ and $q$ are open projections of a $C^*$--algebra such that $\|p-q\|<1$, are the corresponding hereditary subalgebras isomorphic?
It turned out that an affirmative answer with a short proof could be given in the $\sigma$--unital case, which includes the separable case.
It also turned out that special cases of Takesaki's question relate to complete order isomorphisms of $C^*$--algebras, structure of quasi--multipliers, and continuous fields of Hilbert spaces.

The plan of the paper is as follows.
\S2 contains preliminaries and a description of Shen's result.
\S3 contains the basic positive results on Takesaki's question.
\S4 gives an analysis of complete order isomorphisms, a proof that $\sigma$--unital completely order isomorphic $C^*$--algebras are isomorphic, and some results on the structure of quasi--multipliers.
The main results are that in the $\sigma$--unital case quasi--multipliers can be generated from left and right multipliers by multiplication, but not by addition, even for some very nice $C^*$--algebras.
It had previously been shown by McKennon [32] that $QM(A)$ need not be $LM(A)+RM(A)$ in the non--separable case.
\S5 answers a question of Dixmier and Douady [21] by showing that two separable continuous fields of Hilbert spaces which are isomorphic as continuous fields of Banach spaces are also isomorphic as continuous fields of Hilbert spaces.
\S5 also discusses lifting problems for related maps of spaces of operators.
Let $q_0\colon G\to P_0$ be given by $q_0(T)=T^*T$, where $G$ is the space of invertible elements of $B(H)$ with the strong operator topology and $P_0$ is the space of invertible positive operators with the weak operator topology.
Also let $r\colon E\to B$ where $E=B(H)$ with the strong topology, $B$ is the set of self--adjoint operators with the weak topology, and $r(T)={T+T^*\over 2}$.
Then for $H$ separable and infinite dimensional $q_0$ is well behaved and $r$ is badly behaved for lifting problems.
\S6 has examples showing that the positive results of \S3,4 and 5 for the $\sigma$--unital or separable case fail in general.
We also give Example 6.4 which is similar to an example of B.E.~Johnson [26] on perturbations of $C^*$--algebras and accomplishes a slightly stronger result.
In both 6.4 and [26] the $C^*$--algebras are very nice and are isomorphic, but the isomorphism cannot be taken ``small''.
Remark 7.1 discusses the relation between 6.4 and [26] and the relation between the paper as a whole and the subject of perturbations of $C^*$--algebras.

\subhead \S2.\ Preliminaries and description of N.T.~Shen's result\endsubhead

Let $A$ be a $C^*$--algebra and $A^{**}$ its enveloping $W^*$--algebra.
Then $T\in A^{**}$ is called a \underbar{multiplier} of $A$ ($T\in M(A)$) if $TA,AT\subset A$.
Similarly $T$ is a \underbar{left multiplier} $(T\in LM(A))$ if $TA\subset A,T$ is a \underbar{right multiplier} $(T\in RM(A))$ if $AT\subset A$, and $T$ is a \underbar{quasi--multiplier} $(T\in QM(A))$ if $ATA\subset A$.
If $\pi\colon A\to B(H)$ is a faithful representation, then the extension of $\pi$ to $A^{**}$ maps $M(A), LM(A), RM(A)$, and $QM(A)$ isometrically onto the sets of operators in $B(H)$ which satisfy the appropriate multiplication properties relative to $\pi(A)$ (cf.~[34, Proposition 3.12.3]).
Multipliers, etc.~can also be identified with certain maps on $A$ called centralizers (double, left, right or quasi-).
For example, a left centralizer of $A$ is an (automatically bounded) linear map $L\colon A\to A$ such that $L(ab)=L(a)\cdot b$, $\forall a,b\in A$.
There are four topologies natural to use in connection with $M(A)$, etc., and we will regard these topologies as defined on all of $A^{**}$.
The \underbar{strict topology} is generated by the semi--norms $x\mapsto \|xa\|$ and $x\mapsto \|ax\|$, $a\in A$.
Thus a net $(x_\alpha)$ converges to $x$ strictly if and only if $x_\alpha a\to xa$ and $ax_\alpha\to ax$ in norm, $\forall a\in A$.
Similarly, we have the \underbar{left strict topology}, generated by the semi--norms $\|xa\|$, the \underbar{right strict topology}, generated by $\|ax\|$, and the \underbar{quasi--strict topology} generated by $\|a_1xa_2\|$.
For detailed expositions the reader is referred to [3] and [34, \S3.12].

We will make much use of the theory of $A-B$ imprimitivity bimodules, where $A,B$ are $C^*$--algebras.
An $A-B$ imprimitivity bimodule is a vector space $X$ which is a left $A$, right $B$--module and which has two inner products $\langle ,\rangle_A$ and $\langle ,\rangle_B$ taking values in $A$ and $B$ and generating dense subspaces of $A$ and $B$.  Here
$\langle,\rangle_A$ is linear in the first variable, conjugate linear in the second, and respects the $A$--module action on the first variable, and $\langle,\rangle_B$ is conjugate linear in the first variable, linear in the second, and respects the $B$--module action on the second variable.  It is also required that $\langle x,y\rangle_Az=x\langle y,z\rangle_B$.  Let 
$||x||=\|\langle x,x\rangle_A\|^{1/2}=\|\langle x,x\rangle_B\|^{1/2}$, which is a semi--norm with respect to which $X$ is required to be Hausdorff and complete.
See [38], for example, for more details.  
There is a related and weaker concept, $A-B$ Hilbert bimodule, which we will not require in this paper, and also a concept, right Hilbert $B$--module, which we will use a little.
A right Hilbert $B$--module $X$ is endowed only with a right $B$--action and $B$--valued inner product.
It is then possible to define the algebra $A=\sK(X)$ so that $X$ becomes an $A-B$ Hilbert bimodule, and even an $A-B$ imprimitivity bimodule if $\langle X,X\rangle_B$ generates a dense subspace of $B$.
The isomorphism classes of $A-B$ imprimitivity bimodules can be regarded as the morphisms (called strong Morita equivalences) from $A$ to $B$ of a category whose objects are $C^*$--algebras.
All morphisms of this category are invertible, and every ordinary isomorphism induces an imprimitivity bimodule.
Two isomorphisms induce isomorphic imprimitivity bimodules if and only if they differ by an inner automorphism (an automorphism induced by a unitary multiplier).
This point of view (which is used only slightly in this paper) is explained in [8].
If $X$ is an $A-B$ imprimitivity bimodule and $Y$ a $B-C$ imprimitivity bimodule, their composition in the category is a completed tensor product $X\bigotimes_B Y$.

It is often useful to look at an imprimitivity bimodule as a subspace of a $C^*$--algebra in such a way that the four multiplications of the bimodule agree with the algebra multiplication.
Then $\langle x,y\rangle_A$ becomes $xy^*$ and $\langle x,y\rangle_B$ becomes $x^*y$, and 
$X^*=\{x^*\colon x\in X\}$ becomes the inverse $B-A$ imprimitivity bimodule.
Whenever $A$ and $B$ are hereditary $C^*$--subalgebras of a $C^*$--algebra $C$, $X=(ACB)^-$ becomes an $A-B$ Hilbert bimodule in this way.
It is an $A-B$ imprimitivity bimodule if and only if $A$ and $B$ generate the same closed two--sided ideal of $C$.
Every $A-B$ imprimitivity bimodule arises in this way.
In fact, given $A,B,X$ there is a $C^*$--algebra $L$, called the \underbar{linking algebra}, and a projection $p\in M(L)$ such that $A$ is identified with $pLp$, $B$ with $(1-p)L(1-p)$, and $X$ with $pL(1-p)$ ([8]).
Using $L$, one can define $M(X), LM(X)$, etc.~by $M(X)=M(L)\cap pL^{**}(1-p)$, $LM(X)=LM(L)\cap pL^{**}(1-p)$, etc.
Of course $pL^{**}(1-p)$ can be identified with $X^{**}$, and in fact it is possible to make the definitions of $M(X)$, etc.~without mentioning $L$.
The basic ideas of multipliers, etc., including their identification with various kinds of centralizers, all have counterparts in this context.

If $X$ is an $A-B$ imprimitivity bimodule and $Y$ a $B-C$ imprimitivity bimodule, then, as mentioned above, $Z=X\otimes_B Y$ is an $A-C$ imprimitivity bimodule.
It is useful to view this another way.
If $L$ is the linking algebra for $X$, then since $L$ is strongly Morita equivalent to $B,Y$ induces an $L-C$ imprimitivity bimodule $\tilde Y$.
If $\tilde L$ is the linking algebra for $\tilde Y$, we will call $\tilde L$ the \underbar{double linking algebra}.
There are three orthogonal projections $p,q,r\in M(\tilde L)$ such that $p+q+r=1$, $A$ is identified with $p\tilde Lp$, $B$ with $q\tilde Lq$, $C$ with $r\tilde Lr$, $X$ with $p\tilde Lq$, $Y$ with $q\tilde Lr$, and $Z$ with $p\tilde Lr$.
(Also $L$ is identified with $(p+q)\tilde L(p+q)$ and the linking algebra for $Y$ with $(q+r)\tilde L(q+r)$.)
In $\tilde L$ the tensor multiplication for $Z=X\bigotimes_B Y$ becomes algebra multiplication.
It follows that if $T_1\in X^{**}$ and $T_2\in Y^{**}$ we may write $T_1 T_2$ for an element of $Z^{**}\subset \tilde L^{**}$ (instead of $T_1\otimes T_2$).
We are mainly interested in the case where $T_1$ and $T_2$ are quasi--multipliers, at least.

Finally, we mention some notations and conventions that will be used throughout the paper.
The symbol ``$*$'' has different meanings in ``$A^{**}$'' or ``$X^{**}$'' and in ``$X^*$''.
In the one case we have the Banach space double dual, and in the other the adjoint operation.
This should not cause confusion.
It is standard to regard $A^{**}$ as a concrete algebra of operators, namely the double commutant of $A$ in its universal representation.
We will adopt a similar convention with regard to $X^{**}$.
Let the linking algebra $L$ and $p\in M(L)$ be as above.
If $H$ is the Hilbert space for the universal representation of $L$ and $P$ the image of $p$, $X^{**}$ is naturally regarded as a space of operators from $(1-P)H$ to $PH$.
For $T\in X^{**}$, in particular for $T\in QM(X)$, we will take this point of view whenever we use the terms one--one, dense range, invertible, unitary, or isometry.
For example, a unitary element of $X^{**}$ is one which induces an isometry from $(1-P)H$ onto $PH$.
Let $\sK$ denote the $C^*$--algebra of compact operators on a separable infinite dimensional Hilbert space, $\sK(H)$ the compacts on an arbitrary Hilbert space $H$, $\sK(X)$ the compacts of a right Hilbert module $X$, and $\sK(\sH)$ the field of elementary $C^*$--algebras associated with a continuous field $\sH$ of Hilbert spaces ([21]).
A $C^*$--algebra is called $\sigma$--{\it unital} if it posses a strictly positive element, or equivalently if it has a countable approximate identity.
Every separable $C^*$--algebra is $\sigma$--unital.

We now describe N.T.~Shen's result.
Suppose $A$ and $B$ are hereditary $C^*$--subalgebras of $C$.
One wishes to describe the relative position of $A$ and $B$.
A reasonable intuitive concept of what this should mean is that the relative position is already seen in the hereditary $C^*$--subalgebra of $C$ generated by $A\cup B$, and therefore we assume this subalgebra is all of $C$.
Then the $A-B$ Hilbert bimodule $X=(ACB)^-$ is clearly an invariant for the problem.
We will need only the special case where $X$ is an $A-B$ imprimitivity bimodule, though Shen's result is the same for the general case.
One assumes $X$ is given and seeks the other invariants describing the relative position of $A$ and $B$.
Thus we consider embeddings of $A, B$ and $X$ into a $C^*$--algebra $C$, satisfying the above hypotheses and the obvious compatibility conditions.
Two embeddings $\theta\colon (A,B,X)\to C$ and $\theta'(A,B,X)\to C'$ are equivalent if there is an isomorphism of $C$ onto $C'$ which carries $\theta$ onto $\theta'$.  Equivalence of $\theta$ and $\theta'$ is one interpretation of what it means for $\theta(A)$, $\theta(B)$ to have the same relative position as $\theta'(A)$, $\theta'(B)$.

{\bf 2.1.\ Theorem} (Shen [39], also [9, Theorem 4.3).
The equivalence classes of embeddings of $(A,B,X)$ are in one--one correspondence with elements $T$ of $QM(X)$ such that $\|T\|\leq 1$.

The way the correspondence occurs is that for $a\in A$ and $b\in B$, $ \theta(a)\cdot\theta(b)=\theta(aTb)$.
Another way to look at this is that the ways of defining a multiplication on the vector space $A\oplus B\oplus X\oplus X^*$ which lead to an embedding with $C$ the completion of this vector space are parametrized by the contractions in $QM(X)$.
If $p$ and $q$ are the open projections in $C^{**}$ corresponding to $A$ and $B$, then $T$ can be identified with $pq$.
More precisely, the embedding of $X$ into $C$ gives an embedding of $X^{**}$ into $C^{**}$, and the image of $T$ is $pq$.
Under 2.1 the linking algebra corresponds to $T=0$.
It is possible to use other natural equivalence relations to describe the relative position of two hereditary $C^*$--subalgebras and to deduce from 2.1 a description of the equivalence classes.
Shen's result can be regarded as an analogue of the description of the relative position of two subspaces of a Hilbert space ([19]).

Now we use 2.1 to translate Takesaki's question into operator theory.
If $\|p-q\|< 1$, then clearly $A$ and $B$ generate the same ideal of $C$.
Hence $(ACB)^-$ is an $A-B$ imprimitivity bimodule.
Using [19], one sees easily that $\|p-q\|<1$ if and only if $pq$ is invertible as an element of $((ACB)^-)^{**}$.
Thus Takesaki's question becomes:
\item{(T1)}If there is an $A-B$ imprimitivity bimodule $X$ and an invertible $T\in QM(X)$ (such that $\|T\|\leq 1$), does it follow that $A$ is isomorphic to $B$?

The portion of (T1) enclosed in parentheses is clearly of no importance.

We consider also the stronger question:

\item{(T2)}If $X$ is an $A-B$ imprimitivity bimodule such that there is an invertible $T\in QM(X)$, does it follow that $X$ is induced from an isomorphism of $A$ and $B$?

It turns out, not surprisingly, that both questions have the same general answer (yes in the $\sigma$--unital case, no in general), though certainly there are particular examples of $A,B,X$ (not $\sigma$--unital) for which the answer to (T1) is yes and the answer to (T2) is no.

We now give some elaborations of 2.1, of which only the first will be used in this paper.

\medskip
{\bf 2.2.\ Proposition} (cf. [9, Proposition 5.3]).
With the above notations $p\in M(C)$ ($A$ is corner of $C$) if and only if $T\in LM(X)$.
Also $q\in M(C)$ if and only if $T\in RM(X)$.

\medskip
\example{2.3}  (cf.~[9, Theorem 5.4])
Under an embedding $\theta\colon (A,B,X)\to C$ it is possible that $\theta(A)$ and $\theta(B)$ have non--trivial intersection.
We explain how to calculate this intersection from $T$.
Consider the polar decomposition $T=U|T|=|T^*|U$, with $U$ a partial isometry in $X^{**}$, $|T|\in B^{**}$, and $|T^*|\in A^{**}$.
Let $p_1$ and $q_1$ be the spectral projections of $|T^*|$ and $|T|$ for $\{1\}$, $p_2, q_2$ the largest open projections smaller than $p_1,q_1$, and $A', B'$ the corresponding hereditary $C^*$--subalgebras of $A,B$.
Then $\theta(A)\cap\theta(B)=\theta(A')=\theta(B')$.
More precisely, for $a\in A'$, $\theta(a)=\theta(TaT^*)=\theta(UaU^*)$ and $TaT^*=UaU^*\in B'$.
Also $p_2 U=Uq_2$.
\endexample

\example{2.4}
In [1], Akemann proved that if $p$ and $q$ are open projections of $C$ with positive angle, then $p\wedge q$ is open.
Then 2.1 and 2.3 make possible an operator--theoretic translation:\ If $T\in QM(X)$, $\|T\|\leq 1$, and the spectrum of $|T|$ omits an interval $(1-\epsilon,1)$, then the spectral projection of $|T|$ for \{1\} is open.
In connection with this we note that $p\wedge q$ is identified with $p_1$ and $q_1$ in the notation of 2.3.
It is possible to deduce the operator--theoretic translation of Akemann's result from known results on semicontinuous operators.
Nevertheless, we consider that the above is evidence of the usefulness of using 2.1 to translate problems about two hereditary $C^*$--subalgebras into problems of operator theory.
\endexample

Finally, we prove a result generalizing part of a result of Effros, [22, Theorem 2.4], (cf.~also [34, Theorem 1.5.2]).

{\bf 2.5.\ Theorem}.
Let $X$ be an $A-B$ imprimitivity bimodule.
Then there is a one--one correspondence between closed $B$--submodules of $X$ and hereditary $C^*$--subalgebras of $A$, as follows:\ For $X_0$ a closed $B$--submodule the corresponding hereditary subalgebra is $\alpha(X_0)$, the closed linear span of $\langle X_0,X_0\rangle_A$.
For $A_0$ a hereditary $C^*$--subalgebra of $A$ the corresponding submodule is $\sS(A_0)=\overline{A_0 X}$.

\demo{Proof}
Let $A_0 X$ denote the linear span of $\{ax\colon a\in A_0, x\in X\}$.
Clearly, $\sS(A_0)$ is a closed $B$--submodule.
Also $\alpha(\sS(A_0))$ is the closed linear span of $A_0\langle X,X\rangle_A A_0=A_0$ since $A_0$ is hereditary and $\langle X,X\rangle_A$ spans a dense ideal of $A$.
To see that $\alpha(X_0)$ is an algebra, use $\langle x,y\rangle_A \langle z,w\rangle_A=\langle\langle x,y\rangle_A z,w\rangle_A=\langle x\langle y,z\rangle_B,w\rangle_A\in \alpha(X_0)$, $x,y,z,w\in X_0$, since $X_0 B\subset X_0$.
To see that $\alpha(X_0)$ is hereditary, use $\langle x,y\rangle_A a\langle z,w\rangle_A=\langle\langle x,y\rangle_A az,w\rangle_A=\langle x\langle y,az\rangle_B,w\rangle_A\in \alpha(X_0)$.
Now $\sS(\alpha(X_0)$ is the closed span of $\langle X_0,X_0\rangle_A X=X_0\langle X_0,X\rangle_B$.
Clearly, this is a closed $B$--submodule contained in $X_0$.
To show that it is all of $X_0$, we show that $x$ is in the closed $B$--submodule generated by $x\langle x,x\rangle_B$, $x\in X_0$.
This follows from $\langle x,x\rangle_B=|x|^2$.
In any $C^*$--algebra there is a sequence $(p_n)$ of polynomials with no constant term such that $xp_n(|x|^2)\to x$.
\enddemo

\subhead \S3.\ Close open projections:\ positive results\endsubhead

{\bf 3.1.\ Theorem}.
Let $X$ be an $A-B$ imprimitivity bimodule for $A,B$ $\sigma$--unital.
If there is $T\in QM(X)$ which is one--one and has dense range, then $X$ arises from an isomorphism of $A$ and $B$.

\demo{Proof}
Let $e$ and $f$ be strictly positive elements of $A$ and $B$.
Then $x=eTf\in X$ and $x$ is still one--one with dense range.
So $x$ has a polar decomposition, $x=u|x|$, with $u$ a unitary in $X^{**}$.
We claim that $u\in M(X)$.
To see this, note that $x$ one--one implies $|x|$ one--one, which implies that $|x|$ is strictly positive in $B$.
Thus $|x|B$ is dense in $B$, and since $u|x|B=xB\subset X$, it follows that $uB\subset X$.
This shows that $u$ is in $LM(X)$.
That $u$ is in $RM(X)$ can be proved by a similar argument.
By Lemma 3.3 of [8] the existence of a unitary in $M(X)$ is equivalent to the existence of an isomorphism of $A$ and $B$ which induces $X$.
\enddemo

{\bf 3.2.\ Corollary}.
Let $A$ and $B$ be $\sigma$--unital hereditary $C^*$--subalgebras of $C$ and $p$ and $q$ the corresponding open projections.
If $\|p-q\|<1$, then $A$ is isomorphic to $B$.

\medskip
\demo{Proof}
Let $X=(ACB)^-$.
As indicated in \S2, N.T.~Shen's thesis [39] produces a $T\in QM(X)$, and $\|p-q\|<1$ implies that $X$ is an $A-B$ imprimitivity bimodule and $T$ is invertible.
(Recall that $T$ can be identified with $pq$.)
Hence 3.1 applies.
\enddemo

A result of Cuntz (1.4 of [18]) can be derived as a corollary of 3.1.
Cuntz's result could actually be used to prove 3.2, without mentioning imprimitivity bimodules, but we think imprimitivity bimodules provide the right framework for this subject.

\medskip
{\bf 3.3.\ Corollary} (Cuntz [18]).
If $C$ is a $C^*$--algebra, $c\in C$, and $A$ and $B$ are the hereditary $C^*$--subalgebras generated by $cc^*$ and $c^*c$, respectively, then $A$ is isomorphic to $B$.

\demo{Proof}Here $X=(ACB)^- =(cCc)^-$ and $T=c\in X$.
\enddemo

It is of some interest to know whether it is necessary for both $A$ and $B$ to be $\sigma$--unital.
Results showing that sometimes only one need be $\sigma$--unital are presented below.
Relevant counterexamples are 6.1, due to Choi and Christensen [12], which shows that 3.2 fails if neither $A$ nor $B$ is assumed $\sigma$--unital, 6.2, which accomplishes the same purpose for $A,B$ of continuous trace, and 6.3, which shows that 3.1 fails if only one of $A,B$ is assumed $\sigma$--unital.  The proof that 6.3 yields the hypotheses of 3.1 depends on a universal measurability argument such as the one suggested in the paragraph preceding Theorem 4.15 below. 

\medskip
{\bf 3.4.\ Lemma}.
If $\Delta$ is a second countable topological space and $D$ is an upward directed set of lower semicontinuous functions on $\Delta$, then there is a countable $D_0\subset D$ such that $\sup\limits_{g\in D_0} g(t)=\sup\limits_{g\in D} g(t)$, $\forall \ t\in \Delta$.

This is a known result, stated only for reference.

Note that 3.4 applies if $\Delta$ is the quasi-state space of a separable $C^*$-algebra $B$, since the elements of $B$ separate points of $\Delta$ via the Kadison function representation, [34, \S3.10].  Thus the Stone-Weierstrass theorem implies that $C(\Delta)$ is separable, whence $\Delta$ is second countable.

\medskip
{\bf 3.5.\ Theorem}.
Let $X$ be an $A-B$ imprimitivity bimodule with $B$ separable.
If there is $T\in QM(X)$ which  has dense range, then $A$ is $\sigma$--unital.  Thus if also $T$ is one-to-one, then 3.1 applies.

\demo{Proof}
Let $(e_\alpha)_{\alpha\in D}$ be an approximate identity for $A$ and $f$ a strictly positive element for $B$.
For each $\alpha\in D$, $b_\alpha=fT^* e_\alpha Tf\in B$.
By Kadison's function representation each $b_\alpha$ may be regarded as a continuous function on $\Delta$, the quasi--state space of $B$.
The supremum of these functions corresponds to $b=fT^* Tf\in B^{**}$.
By 3.4, we can find $e_{\alpha_1}\leq e_{\alpha_2}\leq\ldots$ such that $b$ is the weak limit of $fT^*e_{\alpha_n}Tf$.
Let $A_0$ be the ($\sigma$--unital) hereditary $C^*$--subalgebra of $A$ generated by the $e_{\alpha_n}$'s and $p$ the corresponding open projection.
Then $fT^*pTf\geq b$.
Hence $(1-p)Tf=0$.
Since $Tf$ has dense range, this shows that $p=1$ and $A=A_0$.
\enddemo

\example{3.6. Remark}
It is trivial that if the hypothesis on $B$ is weakened to $\sigma$--unitality and the hypothesis on $T$ is strengthened to $T\in LM(X)$, then the conclusion still holds.
($TfT^*$ will be strictly positive in $A$.)
In Example 6.3 $T\in RM(X)$.
\endexample

{\bf 3.7.\ Theorem}.
Let $X$ be an $A-B$ imprimitivity bimodule with $B$ $\sigma$--unital.
If there is an invertible $T\in QM(X)$, then $A$ is $\sigma$--unital (and hence 3.1 applies).

\demo{Proof} Let $e_\alpha, D, f, b_\alpha, b$, and $\Delta$ be as in the proof of 3.5.
Let $T^*T\geq\epsilon > 0$, and let $\Delta_n=\{\varphi\in\Delta\colon \varphi (f^2)\geq {1\over n}\}$.
Then $\Delta_n$ is compact and $\bigcup_1^\infty\Delta_n=\Delta\backslash \{0\}$.
Since $b\geq \epsilon f^2$, Dini's theorem (or a similar argument) shows that there is $\alpha_n\in D$ with $b_{\alpha_n}> {\epsilon\over 2} f^2$ on $\Delta_n$ (here $b_\alpha$ and $f^2$ are regarded as functions on $\Delta$).
Let $A_0$ be the ($\sigma$--unital) hereditary $C^*$--subalgebra of $A$ generated by the $e_{\alpha_n}$'s and $p$ the corresponding open projection.
Then $fT^*pTf\geq {\epsilon\over 2} f^2$.
Hence $T^*pT\geq {\epsilon\over 2}$.
Since $T$ is invertible, this implies $p=1$ and $A_0=A$.
\enddemo

\subhead \S4.\ Complete order isomorphisms and the structure of quasi--multipliers\endsubhead

\medskip
{\bf 4.1.\ Lemma}.
If $X$ is an $A-B$ imprimitivity bimodule and $T\in LM(X)$ is invertible, then $T^{-1}\in LM(X^*)$.

\demo{Proof}
$X^*$ is a $B-A$ imprimitivity bimodule.
We need to show that $T^{-1}A\subset X^*$.
Since $T\in LM(X)$, $TB\subset X$, and this implies $TBT^*\subset A$ (using the fact that $B=B^2$).
Clearly $T^{-1}(TBT^*)=BT^*=(TB)^*\subset X^*$.
Thus it is sufficient to show $(TBT^*)^-=A$.
Now $(TBT^*)^-$ is a hereditary $C^*$--subalgebra of $A$ by 2.5.
Also, a state of $A$ vanishes on $TB_+T^*$ if and only if its extension to $A^{**}$ vanishes on $TT^*$, since an approximate identity for $B$ converges to a projection $q$ such that $Tq=T$.
Therefore $(TBT^*)^-=A$ (actually $TBT^*$ is closed since $T$ is invertible), and the proof is complete.
\enddemo

With the hypotheses of 4.1 the map $b\mapsto TbT^*$ is a complete order isomorphism of $B$ onto $A$.
Its inverse is $a \mapsto T^{-1} a(T^{-1})^*$.
Note also that by 2.5 and the proof of 4.1 $(TB)^-=X$.
Since $T$ is invertible, $TB$ is closed and $X=TB$.
This makes possible a correspondence of this situation with something involving only $B$.
Now $h=T^*T$ is an invertible positive element of $QM(B)$.
If $h$ is known, the structure of $X$ as a right Hilbert $B$--module is determined.
Then $A$ is isomorphic to the compacts of this Hilbert $B$--module.  (See [38, Definition 6.4], where it is called the imprimitivity algebra, for a definition of the compacts of a Hilbert module.)  
The structure of $X$ as an $A-B$ bimodule depends on more than just $h$, since the choice of the isomorphism between $A$ and the compacts matters.
To construct a right Hilbert $B$--module from $h$, we let $X_h=B$ as a right $B$--module and define $\langle b_1,b_2\rangle=b_1^* hb_2\in B$.
Clearly if $T$ is given $X_h$ is isomorphic to $X$ as a right Hilbert $B$--module.
Conversely, if only $h$ is given, $X_h$ is a right Hilbert $B$--module.
($\|b\|_{X_h}=\|b^*hb\|^{1/2}\geq \epsilon^{1/2} \|b\|_B$ if $h\geq\epsilon > 0$, and hence $X_h$ is complete.)
The identity map, regarded as a map from $B$ to $X_h$, is a left centralizer of $X_h$; and hence there is $T_h\in LM(X_h)$ such that $X_h=T_h B$ and $T_h^* T_h=h\in QM(B)$.
The fact that $T_h^* T_h$ is invertible implies that $T_h$ is left invertible, and (cf.~2.5) it is then easy to see that $T_h$ is invertible.
Thus there is a correspondence between the situation of 4.1 and invertible positive elements of $QM(B)$.

\medskip
{\bf 4.2.\ Theorem}.

(i)\ Every complete order isomorphism $\varphi$ from $B$ to $A$ arises from a pair $(X,T)$, where $X$ is an $A-B$ imprimitivity bimodule, $T$ is invertible in $LM(X)$, and $\varphi(b)=TbT^*$.

(ii)\ \footnote{Corollary 3.2 of [11] (combined with the results of Kadison) implies (ii) for unital maps and also the fact that a 2--positive anti--isomorphism is an isomorphism.  The definition of 2--positive is given on [11, page 565].  We believe that (ii) is already known in the non--unital case also and that we learned it from a lecture by Choi.
However, we have not found a reference.}
Every 2--order isomorphism is a complete order isomorphism.

(iii)\ $(X_1,T_1)$ and $(X_2,T_2)$ induce the same complete order isomorphism if and only if there is an isomorphism from $X_1$ onto $X_2$ (as imprimitivity bimodules) which carries $T_1$ onto $T_2$.

\demo{Proof}
(i) and (ii):\ Let $\varphi\colon B\to A$ be a 2--order isomorphism.
Then $\varphi^{**}\colon B^{**}\to A^{**}$ is also a 2--order isomorphism and is $\sigma$--weakly continuous.
Let $k=\varphi^{**}(1)$.
Then $k$ is invertible, since $\varphi^{**}$ preserves order units and $A^{**}$ is unital.
The map $\psi\colon A^{**}\to A^{**}$ defined by $\psi(x)=k^{-1/2} xk^{-1/2}$ is a $\sigma$--weakly continuous complete order isomorphism of $A^{**}$ and $\psi\circ \varphi^{**}\colon B^{**}\to A^{**}$ is unital.
By results of Kadison [27], [28] $\psi\circ\varphi^{**}$ is the direct sum of an isomorphism and an anti--isomorphism.
Since an anti--isomorphism is not 2--positive
 (except in the commutative case, when it is also an isomorphism), $\psi\circ\varphi^{**}$ is an isomorphism.
Thus $\psi(\varphi(B))\cdot\psi(\varphi(B))\subset \psi(\varphi(B))$; i.e.~$k^{-1/2} Ak^{-1} Ak^{-1/2}\subset k^{-1/2} Ak^{-1/2}$.
It follows that $Ak^{-1} A\subset A$ and $k^{-1}\in QM(A)$.
As above, with $A$ in place of $B$ and setting $B'=\sK(X_{k^{-1}})$, there is a $B'-A$ imprimitivity bimodule $X_{k^{-1}}$ and an invertible $T_{k^{-1}}\in LM(X_{k^{-1}})$ such that $T_{k^{-1}}{}^*T_{k^{-1}}=k^{-1}$.
Let $\theta\colon A-B'$ be the complete order isomorphism induced by $(X_{k^{-1}},T_{k^{-1}})$.
Then $\theta^{**}(k)=T_{k^{-1}}kT_{k^{-1}}{}^*=T_{k^{-1}} (T_{k^{-1}}{}^* T_{k^{-1}})^{-1} T_{k^{-1}}{}^*=1$.
Thus $\theta\circ\varphi\colon B\to B'$ is a 2--order isomorphism, and $(\theta\circ\varphi)^{**}$ is unital.
Arguing as above from [27], [28], we see that $\theta\circ\varphi$ is a *--homomorphism and hence a *--isomorphism.
We use this isomorphism to identify $B$ with $B'$.
Then $(X_{k^{-1}},T_{k^{-1}})$ gives rise to $(Y,S)$, where $Y$ is a $B-A$ imprimitivity bimodule, $S$ is invertible in $LM(Y)$, and $(Y,S)$ induces $\varphi^{-1}$.
It follows from 4.1 and the discussion after that $\varphi$ is induced by $(X,T)$ where $X=Y^*$ and $T=S^{-1}$.

\NoBlackBoxes
(iii)\ If $(X_1,T_1)$ and $(X_2,T_2)$ induce the same complete order isomorphism $\varphi$, we look at $(\varphi^{**})^{-1}(1)=T_1^{-1} T_1^{-1}{}^*=T_2^{-1} T_2^{-1}{} ^*$.
Since $T_i^{-1} T_i^{-1}{}^*=(T_i{}^*T_i)^{-1}$, we can write $h=T_1{}^*T_1=T_2{}^*T_2$.
Then by the discussion above, there are isomorphisms $\alpha_1,\alpha_2\colon X_1,X_2\to X_h$, as a right Hilbert $B$--modules, such that $\alpha_i$ carries $T_i$ to $T_h$.
It is also clear that there are *--isomorphisms $\beta_1,\beta_2\colon \sK(X_h)\to A$ such that $\varphi(b)=\beta_i(T_h b T_h^*)$.
This implies that $\beta_1=\beta_2$, and the conclusion is clear.
\enddemo

\example{4.3. Remarks}

(i)\ It is not difficult to see what composition looks like in terms of 4.2.
If $(X_1,T_1)$ induces $\varphi_{T_1}\colon B\to A$ and $(X_2,T_2)$ induces $\varphi_{T_2}\colon C\to B$, then $(X_1\otimes_B X_2,T_1T_2)$ induces $\varphi_{T_1}\circ\varphi_{T_2}\colon C\to A$.
(See \S2 for the notation.)

(ii)\ The proof above that $k^{-1}\in QM(A)$ appears to depend on the 2--positivity of $\varphi$, but in fact it would work if $\varphi$ were only an order isomorphism.
In that case $\psi\circ\varphi$ would be a Jordan homomorphism, and we could conclude $ak^{-1}a\in A$, $\forall a\in A$.
If $(e_\lambda)$ is an approximate identity for $A$, then from $e_\lambda k^{-1} e_\lambda\in A$ we deduce $a_1 k^{-1} a_2=\lim (a_1 e_\lambda k^{-1} e_\lambda a_2)\in A$.
Thus the proof of 4.2 shows that every order isomorphism of $C^*$--algebras is the composite of a complete order isomorphism and a Jordan isomorphism.
We suggest that those interested in Jordan algebras should look for the Jordan analogue of imprimitivity bimodules.
\endexample

{\bf 4.4\ Proposition}.
If $(X,T)$ induces a complete order isomorphism from $B$ to $A$ as above, then $X$ arises from an isomorphism of $B$ and $A$ if and only if there is an invertible $S\in LM(B)$ such that $S^*S=T^*T$.

\demo{Proof} iF $X$ arises from an isomorphism, then by [8] there is a unitary $U\in M(X)$.
Then $S=U^*T$ has the desired property.

Conversely, if $S$ is given, then $S^{-1}\in LM(B)$ by 4.1.
Then $U=TS^{-1}$ will be a unitary in $LM(X)$.
Every unitary left multiplier is a multiplier.
This follows from Proposition 4.4 of [2] or can be deduced from 4.1.
Thus, $U\in M(X)$ and [8] implies that $X$ arises from an isomorphism.
\enddemo

{\bf 4.5.\ Definition}.
A complete order automorphism, $\varphi\colon B\to B$, will be called \underbar{inner} if there is invertible $S\in LM(B)$ such that $\varphi(b)=SbS^*$.

We now state three properties which a $C^*$--algebra $B$ might satisfy:

\noindent
(P1)\ Every $C^*$--algebra completely order isomorphic to $B$ is isomorphic to $B$.

\noindent
(P2)\ Every invertible positive quasi--multiplier of $B$ is of the form $S^*S$ for an invertible $S$ in $LM(B)$.

\noindent
(P2$'$)\ Every complete order isomorphism of $B$ with a $C^*$--algebra is the composite of an inner complete order automorphism of $B$ and an isomorphism.

It is clear from the above and the discussion in \S2 that the questions whether $B$ satisfies (P1) and (P2) are special cases of (T1) and (T2).
The special case in question is, by 2.2, the case where the open projection $p$ is actually a multiplier (i.e., $A$ is a corner of $C$).
It is interesting that in this case the complete order isomorphism can easily be described in terms of $C\colon\varphi (b)=pbp\in A$, $\forall b\in B$.
Obviously (P2) $\Rightarrow$ (P1).
From the proof of 4.4 it should be fairly clear that (P2) $\Leftrightarrow$ (P2$'$).
One only has to note that $T=US\Rightarrow \varphi_T=\varphi_U\circ\varphi_S$; and that the complete order isomorphism induced by $(X,T)$ is an isomorphism if and only if $T$ is unitary.
It is somewhat amusing to note that, by (P2) $\Rightarrow$ (P1), even for non--$\sigma$--unital $C^*$--algebras property (P2) is preserved by complete order isomorphisms.

The counterexamples 6.1 and 6.2 are actually counterexamples to (P1) as well as (T1).
The positive results of \S3 now yield corollaries:

\medskip
{\bf 4.6.\ Corollary}.
Every $\sigma$--unital $C^*$--algebra satisfies (P2$'$) and (P1).

\medskip
{\bf 4.7.\ Corollary}.
If $B$ is a $\sigma$--unital $C^*$--algebra, then any $C^*$--algebra order isomorphic to $B$ is Jordan isomorphic to $B$.

4.7 is analogous to 4.6 and relies on 4.3(ii).

\medskip
{\bf 4.8.\ Corollary}.
Any $\sigma$--unital $C^*$--algebra satisfies (P2).

4.8 is a result on the structure of quasi--multipliers and is the starting point for the rest of this section.

\medskip
{\bf 4.9.\ Theorem}.
If $A$ is $\sigma$--unital and $0\leq h\in QM(A)$, then there is $T\in LM(A)$ such that $T^*T=h$.

\demo{Proof} Let $e$ be strictly positive in $A$.
Let $a=(ehe)^{1/2}\in A$.
Since $a^*a\leq \|h\|e^2$, it follows that there are $U_n\in A$ with $\|U_n\|\leq \|h\|^{1/2}$ and $U_n e\to a$ in norm $(U_n=a(e^2+{1\over n})^{-1}e$, cf.~[34, p.~12]).
Then $(U_n x)$ is norm convergent for $x\in eA$.
Since $\|U_n\|$ is bounded and $\overline{eA}=A$, it follows that $U_n$ converges in the left strict topology to some $T\in LM(A)$ and $Te=a$.
Therefore $eT^*Te=a^*a=ehe$, and $T^*T=h$.
\enddemo

Note that if the $h$ in 4.9 is one--one, then the $T$ produced in the proof has dense range.
($T$ is one--one also since $T^*T=h$.)
Thus 4.8 could be deduced from 4.9, giving a short proof free of imprimitivity bimodules.
(The author believes that the concept of imprimitivity bimodule provides valuable motivation for this subject and should not be avoided, even when it could be avoided.)

We have already mentioned that the example of Choi and Christensen [12] shows that (P2) may fail if $B$ is not $\sigma$--unital.
In 6.1 we show that this example accomplishes more:\ There is an invertible positive $h\in QM(B)$ such that $h\neq T^*T$, $\forall T\in LM(B)$ (i.e., $T$ is not required to be invertible).

{\bf 4.10.\ Theorem}.
If $A$ is $\sigma$--unital, $R\in QM(A)$, and $\|R\|\leq 1$, then there are isometric columns $L_1=\bmatrix L_{11}\\ L_{21}\endbmatrix$ and $L_2=\bmatrix L_{12}\\ L_{27}\endbmatrix$ of left multipliers such that $R=L_1 {}^*L_2$ (i.e.~$R=L_{11}{}^* L_{12}+L_{21}{}^* L_{27}$, $L_{ij}\in LM(A))$.
If $A$ is stable, then $R=L'_1 {}^* L'_2$ with $L'_1$ and $L'_2$ isometries in $LM(A)$.

\demo{Proof} We first discuss the notation.
Let $X=A\oplus A$, regarded as an $(A\otimes M_2)-A$ imprimitivity bimodule.
$L_1$ and $L_2$ should really be regarded as isometries in $LM(X)$.
If $A$ is stable $(A\simeq A\otimes\sK)$, then there is a unitary $U\in M(X)$, and we deduce the second sentence from the first by writing $L'_i=U^* L_i$.

Now we proceed by applying 4.9 to $h=\bmatrix 1&R\\ R^*&1\endbmatrix\in QM (A\otimes M_2)$.
Let $L=(L_{ij})\in LM(A\otimes M_2)$ be such that $L^* L=h$.
Then if $L_1$ and $L_2$ are the columns of $L$, $L_i{} ^*L_i=1$ and $L_1{} ^*L_2=R$.
\enddemo

\example{4.11. Remarks}

(i)\ In [4], we proved the same result for $R$, $L'_i\in B(H)$.
The present proof is better, we think.

(ii)\ Examples 6.5 and 4.22 below show that $R$ may not be the product of two invertible elements of $RM(A)$ and $LM(A)$, even if $R$ and $R^{-1}\in QM(A)$.
\endexample

Consider $RM(A)\hat\otimes LM(A)$, where the tensor product is given the maximal cross--norm and then completed.
Since $RM(A)\cdot LM(A)\subset QM(A)$ and $QM(A)$ is norm closed, multiplication gives a map $\mu\colon RM(A)\hat\otimes LM(A)\to QM(A)$.

{\bf 4.12.\ Corollary}.
If $A$ is $\sigma$--unital, then $\mu$ is surjective.
If $\overline\mu$ is the induced map from $RM(A)\hat\otimes LM(A)/{ker\ \mu}$ to $QM(A)$, then $\|\overline\mu^{-1}\|\leq 2$.
If $A$ is stable, $\overline\mu$ is an isometry.

Let $\pi\colon A\to B$ be a surjective homomorphism of $C^*$--algebras.
Then $\pi$ induces maps $M(A)\to M(B)$, $LM(A)\to LM(B)$, etc.
In [3] Akemann, Pedersen, and Tomiyama showed that the map $M(A)\to M(B)$ is surjective if $A$ is separable.
Pedersen [35] shows that this is true if $A$ is only $\sigma$--unital.
The corresponding results for $LM$ and $QM$ are proved below.

{\bf 4.13.\ Theorem}.
If $\pi\colon A\to B$ is a surjective $*$--homomorphism of $C^*$--algebras, where $A$ is $\sigma$--unital, and $T\in LM(B)$, then there is $\tilde T\in LM(A)$ such that $\pi^{**}(\tilde T)=T$ and $\|\tilde T\|=\|T\|$.

\demo{Proof} Let $e$ be strictly positive in $A$.
Let $x=T\pi(e)\in B$.
Then $x^*x\leq \|T\|^2 \pi(e^2)$.
By 1.5.10 of [34] (a result of Combes [16]), there is a $y\in A$ with $y^*y\leq \|T\|^2 e^2$ and $\pi(y)=x$.
Then there are $u_n\in A$ such that $\|u_n\|\leq \|T\|$ and $u_n e\to y$, in norm.
As in the proof of 4.9, there is $\tilde T\in LM(A)$ such that $u_n\to\tilde T$ in the left strict topology.
Thus $\|\tilde T\|\leq \|T\|$ and $\tilde T e=y$.
Hence $\pi^{**}(\tilde T)\pi (e)=\pi(y)=T\pi (e)$.
Since $\pi(e)$ is strictly positive in $B$, $\pi^{**}(\tilde T)=T$.
\enddemo

{\bf 4.14.\ Theorem}.
If $\pi,A$ and $B$ are as in 4.13 and $T\in QM(B)$, there is $\tilde T\in QM(A)$ such that $\|\tilde T\|=\|T\|$ and $\pi^{**}(\tilde T)=T$.

\demo{Proof} We first prove the result for $B$ stable.
If $\|T\|=1$, 4.10 shows that $T=L_1{}^* L_2$ with $L_i\in LM(B)$ and $\|L_i\|=1$.
By 4.13, $L_i=\pi^{**}(\tilde L_i)$ with $\tilde L_i\in LM(A)$ and $\|\tilde L_i\|=1$.
Set $\tilde T=\tilde L_1 {}^* \tilde L_2$.

Now in the general case we can consider $\pi\otimes id\colon A\otimes\sK\to B\otimes\sK$.
Let $p$ be a rank one projection in $\sK$ and identify $A$ and $B$ with $A\otimes p$ and $B\otimes p$.
Then there is $\tilde T'\in QM(A\otimes\sK)$ with $(\pi\otimes id)^{**} (\tilde T')=T\otimes p$ and $\|\tilde T'\|=\|T\|$.
Define $\tilde T$ by $\tilde T\otimes p=(1\otimes p)\tilde T'(1\otimes p)$.
\enddemo

We next prove a fairly technical result, that the $C^*$--algebra generated by $QM(A)$ is contained in $\sU(A)+i\sU(A)$, where $\sU(A)$ is the space of universally measurable elements of $A^{**}$ ([33] or [34, 4.3.11]).
Our motivation for this is that the atomic representation of $A$ is faithful on $\sU(A)$ ([34, 4.3.15]), and also, by [34, 4.5.12], for any $T\in A^{**}$ such that the $C^*$--algebra generated by $T$ is in $\sU(A)+i\sU(A)$, the domain and range projections of $T$ are in $\sU(A)$.
Hence one can check whether $T$ is one--one or has dense range by using only irreducible representations of $A$, rather than the universal representation.
This is useful when $A$ is the algebra of continuous sections vanishing at $\infty$ of $\sK(\sH)$ for a continuous field of Hilbert spaces $\sH$, since then $A$ is presented in terms of its universal atomic representation.
We are adopting the notations of [34], and in particular sections 3.11, 4.3 and 4.5 of [34] are prerequisites.
We denote $(A_{sa}^m-A_{sa}^m)^-$ by $\tilde\sB_0$.
Then $\tilde\sB_0$ is a real vector space contained in $\sU(A)$ and $\tilde\sB_0\supset\tilde A_{sa}^m$.
It was shown by Combes [16] that $\tilde\sB_0$ is a Jordan algebra.

{\bf 4.15.\ Theorem}.
If $Q$ is the $C^*$--algebra generated by $QM(A)$, then $Q\subset\tilde\sB_0+i\tilde\sB_0\subset \sU(A)+i\sU(A)$.

\demo{Proof} Let $T\in QM(A)$ and $a\in A_{sa}$.
Then $T^* aT\in QM(A)_{sa}\subset((\tilde A_{sa})^m)^-$, by [34, 3.12.9].
It follows that $h\in A_{sa}^m\Rightarrow T^*hT\in ((\tilde A_{sa})^m)^-$ (it follows from [34, 3.11.8] that $((\tilde A_{sa})^m)^- = (((\tilde A_{sa})^m)^-)^m)$.
Hence $T^*\tilde\sB_0 T\subset\tilde\sB_0$.
By polarization, if $S,T\in QM(A)$ and $x\in\tilde\sB_0$, $T^*x S+ S^*x T\in\tilde\sB_0$.
Choosing $S=1$, we see that $T^*x+xT\in\tilde\sB_0$, and choosing $S=i$, we see that $T^*x-xT\in i\tilde\sB_0$.
Thus $T^*x$, $xT\in\tilde\sB_0+i\tilde\sB_0$.
Since $QM(A)=QM(A)^*$, this shows that $QM(A)(\tilde\sB_0+i\tilde\sB_0)\subset\tilde\sB_0+i\tilde\sB_0$ and $(\tilde\sB_0+i\tilde\sB_0) QM(A)\subset \tilde\sB_0+i\tilde\sB_0$.
Note that $\sI=\{T\in A^{**}\colon T(\tilde\sB_0+i\tilde\sB_0)\subset\tilde\sB_0 +i\tilde\sB_0$ and $(\tilde\sB_0+i\tilde\sB_0)T\subset \tilde\sB_0+i\tilde\sB_0\}$ is a $C^*$--algebra.
Thus $Q\subset\sI$, and since $\tilde\sB_0+i\tilde\sB_0$ contains 1, $Q\subset\tilde\sB_0+i\tilde\sB_0$.
\enddemo

{\bf 4.16.\ Theorem}.
If $h\in QM(A)$ is positive and invertible, then there is a net $(h_\alpha)_{\alpha\in D}$ such that $h_a\nearrow h$, $h_\alpha\to h$ in the strict topology, and each $h_\alpha=T_\alpha{}^* T_\alpha$ for some $T_\alpha\in LM(A)$.
Moreover, each $h_\alpha$ is of the form $ha_\alpha h$, $a_\alpha\in A_+$, and if $A$ is $\sigma$--unital, $(h_\alpha)$ can be taken as a sequence.

\demo{Proof} By the discussion following 4.1 there are a $B-A$ imprimitivity bimodule $X$ and an invertible $S\in LM(X)$ such that $S^*S=h$.
Since $SAS^*=B$, we can choose $a_\alpha\in A$ such that $Sa_\alpha S^*$ is an (increasing) approximate identity of $B$.
Then $ha_\alpha h=S^* Sa_\alpha S^*S\nearrow S^* S=h$.
We take $T_\alpha=a_\alpha^{1/2} h$.
To see that the convergence is strict, let $a\in A$.
Then $Sa\in X$.
Since $X=BX$ and $Sa_\alpha S^*$ is an approximate identity of $B$, $Sa_\alpha S^* Sa\to Sa$ in norm.
Hence $h_\alpha a=S^* Sa_\alpha S^* Sa\to S^*Sa=ha$ in norm.
By symmetry, $ah_\alpha\to ah$ in norm.
For the second sentence, it is necessary only to point out that $A$ $\sigma$--unital $\Rightarrow B$ $\sigma$--unital $\Rightarrow B$ has a countable approximate identity.
\enddemo

\example{4.17. Remark}
The conclusions of the first sentence of 4.16 are of course useless in the $\sigma$--unital case, since then $h$ is already of the form $T^*T$, $T\in LM(A)$.
The reason we have included the second sentence, and mentioned the $\sigma$--unital case in it, is that $ha_\alpha h\nearrow h\Rightarrow a_\alpha\nearrow h^{-1}$.
Thus for $h$ positive and invertible, $h\in QM(A)$ implies $h^{-1}\in A_+ ^m$.
Now the equivalence of (i) and (iii) in [34, 3.11.8] yields that $h\in ((\tilde A_{sa})^m)^-$ implies $h^{-1}\in (A_{sa}^m)^-$.
Thus we have deduced a stronger conclusion from 4.16, with a stronger hypothesis.
It is unknown whether $A_{sa}^m=(A_{sa}^m)^-$.
We are not sure whether this remark is suggestive with regard to this question.
\endexample

{\bf 4.18.\ Corollary}.\ If $LM(A)=M(A)$, then $QM(A)=M(A)$.

\demo{Proof} It is sufficient to prove that $h\in M(A)$ when $h$ is as in 4.16.
If $h_\alpha$ is as in 4.16, then since $T_\alpha\in M(A)$, also $h_\alpha\in M(A)$.
Since $h_\alpha\to h$ strictly, $h$ must be in $M(A)$.
\enddemo

The theme of much of the above has been to show that $QM(A)$ is generated by $LM(A)$ and $RM(A)$ in some sense.
For $A$ $\sigma$--unital 4.9 and 4.10 are strong results of this type.
For $A$ general 4.16 is a rather weak result.
It is trivial that if $e_\alpha$ is an approximate identity of $A$ and $h\in QM(A)$, then $e_\alpha he_\alpha\in A$ and $e_\alpha he_\alpha\to h$ quasi--strictly.
Also $he_\alpha\to h$ left strictly and $he_\alpha\in RM(A)$.
Comparison of these trivial results with 4.16 makes one wonder whether 4.16 accomplishes much with regard to the theme mentioned above.
4.18 shows that 4.16 is good for something.
In the remainder of this section we consider another attempt to generate $QM(A)$ from $LM(A)$ and $RM(A)$.
Namely the question was raised by Akemann and Pedersen [2] whether $QM(A)=LM(A)+ RM(A)$.
McKennon [32] gave a non--separable counterexample.
We show below that the answer is negative in the separable case also.
It is clear \underbar{a priori} that the answer is most likely to be negative for stable $C^*$--algebras, and it turns out that it is quite typical for the answer to be negative in the stable case.

Let $\alpha=(a_{ij})^\infty_{i,j=1}$ be an infinite matrix.
If $n_1<n_2<\ldots,\alpha$ can be subdivided into blocks of size $(n_k-n_{k-1})\times (n_\ell-n_{\ell-1})$ (set $n_0=0$).
Let $L(\alpha)=(b_{ij})$ where 
$$
b_{ij}=\cases a_{ij},&\text{if $\exists k>\ell$ such that $n_{k-1} < i\leq n_k,\ n_{\ell-1} < j\leq n_\ell$}\\
0,&\text{otherwise}. \endcases
$$
Thus $L(\alpha)$ is the lower triangular part of the block matrix.
Similarly, we define $L_0(\alpha)=(c_{ij})$ by
$$
c_{ij}=\cases a_{ij},&\text{if $\exists k>\ell+1$ such that $n_{k-1} < i\leq n_k,\ n_{\ell-1} < j\leq n_\ell$}\\
0,&\text{otherwise}.\endcases
$$
Thus $L(\alpha)-L_0(\alpha)$ consists of only one block--diagonal.
We define $U(\alpha)=L(\alpha^*)^*$.

Let $A$ be a unital $C^*$--algebra and $B=A\otimes\sK$.
If $\{e_{ij}\}^\infty_{i,j=1}$ is a set of matrix units for $\sK$, an element $b$ of $B^{**}$ can be identified with an infinite matrix $(a_{ij})$, where $a_{ij}\in A^{**}$ is defined by $(1\otimes e_{ii}) b(1\otimes e_{jj})=a_{ij}\otimes e_{ij}$.

{\bf 4.19.\ Lemma}.
Let $A$ be a unital $C^*$--algebra, $B=A\otimes\sK$, and $\alpha=(a_{ij})$ an infinite matrix over $A^{**}$.

(i)\ $\alpha$ represents an element of $QM(B)$ if and only if $\alpha$ is bounded and each $a_{ij}\in A$.

(ii)\ $\alpha$ represents an element of $LM(B)$ if and only if $\alpha$ is bounded, each $a_{ij}\in A$, and there exists $n_1< n_2<\ldots$ such that $L_0(\alpha)$ represents an element of $B$.

\demo{Proof}

(i)\ Since $1\otimes e_{ii}$, $1\otimes e_{jj}\in B$, it is obvious that the condition is necessary.
For the converse note that $A\otimes_{alg} \Cal F$ is dense in $B$, where $\Cal F$ is the linear span of $\{e_{ij}\}$.
If $\alpha$ is bounded, to verify that $\alpha$ represents a quasi--multiplier, it is sufficient to check $(A\otimes_{alg}\Cal F) \alpha(A\otimes_{alg}\Cal F)\subset B$.

(ii)\ It follows from (i) that the first two parts of the condition are necessary.
For the third part we assume $\alpha\in LM(B)$.
Write $f_n=\sum_1^n 1\otimes e_{ii}$.
We recursively define $n_1< n_2< \ldots$ so that $\|(1-f_{n_{k+1}})\alpha f_{n_k}\| < 2^{-k}$.
This is possible since $\alpha f_{n_k}\in B$ and $(f_n)^\infty_{n=1}$ is an approximate identity for $B$.
Then $L_0(\alpha)$ is a norm convergent sum, $\sum_1^\infty (1-f_{n_{k+1}})\alpha (f_{n_k}-f_{n_{k-1}})$, and each term of this sum is in $B$.

For the converse, write $\alpha=L_0(\alpha)+\beta$, where $L_0(\alpha)$ is given to be in $B$.
Since by hypothesis $\alpha$ and $L_0(\alpha)$ are bounded, $\beta$ is also bounded.
Since $B\subset LM(B)$, it is sufficient to show $\beta\in LM(B)$; and for this it is sufficient to check $\beta(A\otimes_{alg}\Cal F)\subset B$, which is clear.
\enddemo

For an element $t$ of a $C^*$-algebra, $\text{Re } t$ denotes $(t+t^*)/2$, and for a subset $S$, $\text{Re } S$ denotes $\{\text{Re } t\colon t\in S\}$.

{\bf 4.20.\ Theorem}.
Let $A$ be a unital $C^*$--algebra and $B=A\otimes\sK$.
Then $QM(B)=LM(B)+RM(B)$ if and only if for every bounded self--adjoint matrix $\alpha=(a_{ij})^\infty_{i,j=1}$, with $a_{ij}\in A$, there exists $n_1< n_2<\ldots$ such that $L(\alpha)$ is bounded.

\demo{Proof} Obviously $QM(B)=LM(B)+RM(B)\Leftrightarrow QM(B)_{sa}=\text{Re}\ LM(B)$.
Assume $QM(B)_{sa}=\text{ Re }LM(B)$ and let $\alpha=(a_{ij})$ be a bounded self--adjoint matrix over $A$.
By 4.19(i) $\alpha\in QM(B)$.
Thus there exists a matrix $\beta\in LM(B)$ such that $\alpha=\beta+\beta^*$.
Let $n_1< n_2<\ldots$ be chosen so that $L_0(\beta)\in B$, and write $\beta=L_0(\beta)+\gamma$.
Then $L(\alpha)=L_0(\beta)+L(\gamma)+U(\gamma)^*$.
Then $L_0(\beta)$ is bounded since $L_0(\beta)\in B$, and 
$L(\gamma)$ is bounded since $\gamma$ is bounded and $L(\gamma)$ is a single block--diagonal of $\gamma$.
Also $U(\gamma)$ is bounded, since $\gamma$ is bounded and $\gamma-U(\gamma)$ consists of just two block--diagonals of $\gamma$.

Next assume the condition of the theorem and let $\alpha\in QM(B)_{sa}$.
Then $\alpha$ is a bounded self--adjoint matrix over $A$.
Choose $n_1< n_2 <\ldots$ such that $L(\alpha)$ is bounded.
Let $\sigma=-L(\alpha)+ L(\alpha)^*$ and $\beta=\alpha+\sigma$.
Then $\beta$ is bounded, \text{Re} $\beta=\alpha$, and $L(\beta)=0$.
This implies $L_0(\beta)=0$, so that by 4.19(ii) $\beta\in LM(B)$.
Therefore $QM(B)_{sa}=\text{ Re }LM(B)$.
\enddemo

{\bf 4.21.\ Corollary}.
If $A$ is a unital $C^*$--algebra such that $QM(A\otimes\sK)=LM(A\otimes\sK)+RM(A\otimes\sK)$ and $A_0$ is a unital $C^*$--subalgebra of $A$, then $QM(A_0\otimes\sK)=LM(A_0\otimes\sK)+RM(A_0\otimes\sK)$.

\medskip
{\bf 4.22.\ Example}.
Let $S_1=\{h\in B(\ell^2)\colon h^*=h$ and $\|h\|\leq 1\}$, with the weak operator topology.
Then $S_1$ is a compact metric space.
Let $X$ be any compact Hausdorff space containing a non--empty perfect subspace and let $C$ be the Cantor set.
Then $X$ has a closed subset $X_0$ which can be mapped onto $C$.
Also $C$ can be mapped onto every compact metric space.
Thus $X_0$ can be mapped onto $S_1$.
Since $S_1$ is convex, the Dugundji extension theorem implies $S_1$ is an absolute retract; and hence a map of $X_0$ onto $S_1$ can be extended to $X$.
We now take $A=C(X)$, $B=A\otimes\sK$, and claim that $QM(B)\neq LM(B)+RM(B)$.
Note that $X$ could be a very nice space such as $C$ or $[0,1]$.

To prove the claim we use 4.20.
Lemma 4.19 gives an identification of $QM(B)$ with bounded matrices over $C(X)$.
It is well known that $QM(B)$ can also be identified with the set of weakly continuous functions from $X$ to $B(\ell^2)$ (these are automatically bounded since $X$ is compact).
The relationship between these identifications is the obvious one based on the identification of $B(\ell^2)$ with bounded matrices over $\bC$.
Then let $\alpha\in QM(B)_{sa}$ correspond to a surjective map from $X$ to $S_1$.
If there were $n_1< n_2<\ldots$ such that $L(\alpha)$ is bounded, then for this choice of the $n_k$'s $L(\beta)$ would be bounded for \underbar{every} matrix $\beta$ over $\bC$ representing an element of $S_1$.
Now it is well known (see, for example, [25], where this fact is also used) that there is a matrix $\beta_1=(t_{ij})$, representing an element of $S_1$, such that $L(\beta_1)$, defined relative to $n_k=k$, is not bounded.
Take $\beta=(s_{ij})$ where
$$
s_{ij}=\cases t_{k\ell},&\text{if $i=n_k,j=n_\ell$}\\
0,&\text{otherwise}.\endcases
$$
Then $\beta\in S_1$ and $L(\beta)$ is not bounded.
Thus 4.20 implies $QM(B)\neq LM(B)+RM(B)$.

With the help of 4.21 we can extend the negative result of 4.22.
For simplicity, we state the result only in the separable case.

{\bf 4.23.\ Theorem}.\ If $A$ is a separable unital $C^*$--algebra, then $QM(A\otimes\sK)=LM(A\otimes\sK)+RM(A\otimes\sK)$ implies $A^*$ is separable, or equivalently (by [23, Theorem 3.1]), $A$ has a composition series with elementary quotients.

\demo{Proof} By 4.14 the property $QM(A\otimes\sK)=LM(A\otimes\sK)+RM(A\otimes\sK)$ carries over to quotients of $A$.
Therefore it is sufficient to show $A\neq 0$ implies $A$ has an elementary ideal.
To do this, take $A_0=C(X)$ a MASA in $A$.
Then 4.21 and 4.22 imply $X$ has no perfect subset.
In particular $X$ has an isolated point $x_0$.
It is well known that $x_0$ gives rise to a minimal projection $e_0\in A$ and that the ideal generated by $e_0$ is elementary.
\enddemo

\subhead \S5.\ Continuous fields of Hilbert spaces and related fiber maps\endsubhead

If the results of \S3 and 4 are specialized to the case of (continuous trace) $C^*$--algebras derived from continuous fields of Hilbert spaces, results on continuous fields are obtained.
However, for two reasons such results are not optimal, and one should prove analogues of the results in \S3 and 4 rather than deriving corollaries.
The reasons are:

1.\ Continuous fields can be studied over arbitrary base spaces, not just locally compact spaces.

2.\ The results of \S3 and 4 depend on $\sigma$--unitality.
For a continuous field of Hilbert spaces $\sH$ defined on a locally compact space $X$, the associated $C^*$--algebra is $\sigma$--unital if and only if $\sH$ is separable and $X$ is $\sigma$--compact.
In most cases only the separability of $\sH$ is needed for the results of this section.

In order to avoid devoting too much space to technicalities and redundancies, we have tried to make our proofs in this section brief.

\medskip
{\bf 5.1.\ Definition}.
Let $\sH_1$ and $\sH_2$ be continuous fields of Hilbert spaces on $X$, and let $\sK(\sH_1)$ and $\sK(\sH_2)$ be the associated fields of elementary $C^*$--algebras.
If $T$ is a locally bounded function such that $T(x)\in B(\sH_1(x),\sH_2(x))\ \forall \ x\in X$, then $T$ is called an $L$--map if $e$ a section of $\sH_1\Rightarrow Te$ is a section of $\sH_2$; an $R$--map if $T(\cdot)^*$ is an $L$--map (from $\sH_2$ to $\sH_1$), and a $Q$--map if $e_i$ a section of $\sH_i\Rightarrow (Te_1,e_2)$ is continuous.  Here $Te$ is the function $x\mapsto T(x)e(x)$, $x\in X$.

Then $T$ is both an $L$--map and an $R$--map if and only if $T$ is a map in the usual sense, and 
$T$ is an $L$--map if and only if $T$ is a map in the sense of continuous fields of Banach spaces.
It is possible to prove a precise result showing that $L$--maps, etc.~are the correct analogues of left multipliers, etc., but we content ourselves with the following lemma.

\vfill\eject
\noindent

{\bf 5.2.\ Lemma}.

(i)\ If $T$ is an $L$--map and $k_1$ is a section of $\sK(\sH_1)$, then $Tk_1$ is a section of $\sK(\sH_1,\sH_2)$.

(ii)\ If $T$ is an $R$--map and $k_2$ is a section of $\sK(\sH_2)$ then $k_2T$ is a section of $\sK(\sH_1,\sH_2)$.

(iii)\ If $T$ is a $Q$--map and $k_i$ is a section of $\sK(\sH_i)$, then $k_2 Tk_1$ is a section of $\sK(\sH_1,\sH_2)$.
Also if $f_1$ is a section of $\sH_1, k_2 Tf_1$ is a section of $\sH_2$.

\demo{Proof} Since $T$ is locally bounded, it is sufficient to consider the case $k_i=\sum^n_{j=1} e_j\times f_j$, where $e_j,f_j$ are vector sections.
\enddemo

{\bf 5.3.\ Theorem} (cf.~3.1).\ If $\sH_1$, and $\sH_2$ are separable and $T$ is a $Q$--map such that $T(x)$ is one--one with dense range, $\forall \ x\in X$, then $\sH_1$ and $\sH_2$ are isomorphic.

\demo{Proof} There are sections $e_i$ of $\sK(\sH_i)$ such that $e_i(x)$ is one--one with dense range (cf.~ [16, proof of 10.8.5]).
Thus $e_2 Te_1$ is a section of $\sK(\sH_1,\sH_2)$; and if $e_2(x)T(x)e_1(x)=U(x)h(x)$ is its polar decomposition, $U$ is a unitary map, giving the desired isomorphism.
\enddemo

\example{5.4. Remarks}

(i) An affirmative answer to question $2^\circ$, p.~265 of [21], in the special case of separable fields, follows from 5.3.
Question $2^\circ$ asks about the case where $T$ is an $L$--map with an inverse which is also an $L$--map.
A restatement of our result is that two separable continuous fields of Hilbert spaces which are isomorphic as continuous fields of Banach spaces are isomorphic as continuous fields of Hilbert spaces.
The answer to $2^\circ$ is negative in the non--separable case, as is shown by example 6.2.

(ii)\ The hypothesis that there is a section $e$ of $\sK(\sH)$ such that $e(x)$ is one--one with dense range, $\forall \ x\in X$, is equivalent to separability of $\sH$ if $X$ is paracompact, but is weaker than separability in general.
Thus the separability hypotheses of 5.3 and some of the results below could be weakened.
\endexample

{\bf 5.5\ Theorem} (cf.~3.5).
Assume $\sH_1$ is separable, $T$ is a $Q$--map from $\sH_1$ to $\sH_2$ such that $T(x)$ is one--one with dense range, $\forall \ x\in X$, and $X$ is second countable.
Then $\sH_2$ is separable (and hence 5.3 applies).

\demo{Proof} There is a net $(e_\alpha)_{\alpha\in D}$ of sections of $\sK(\sH_2)$ such that each $e_\alpha=\sum_1^{n_\alpha} g_{j\alpha}\times g_{j\alpha}$, where the $g_{ja}$'s are sections of $\sH_2$, and $0\leq e_\alpha(x)\nearrow 1$, $\forall \ x\in X$.
The proof of this is similar to the proof of the existence of approximate identities in $C^*$--algebras (see [20] or [34]).
Let $f_1,f_2,\ldots$ be a sequence of sections of $\sH_1$, total at each $x\in X$.
Consider $h_{\alpha,k}=(e_\alpha Tf_k,Tf_k)$.
By 3.4 there is a countable $D_0\subset D$ such that $\sup\limits_{\alpha\in D_0} h_{\alpha,k}(x)=(Tf_k(x),Tf_k(x))$, $\forall x,k$.
This shows that the $g_{j\alpha}'s$, $\alpha\in D_0$, are total at each $x\in X$.
\enddemo

\example{5.6. Remark} (cf.~3.6).
If $\sH_1$ is separable and $T$ is an $L$--map from $\sH_1$ to $\sH_2$ such that $T(x)$ has dense range $\forall x$, then $\sH_2$ is separable (no hypothesis on $X$).
Example 6.3 shows that (with $T$ an $R$--map and $X$ compact), the second countability hypothesis cannot be eliminated from 5.5.
\endexample

{\bf 5.7.\ Theorem} (cf.~3.7).
Assume $X$ is paracompact, $T$ is a $Q$--map from $\sH_1$ to $\sH_2$ such that $T(x)^{-1}$ exists, $\forall \ x\in X,\ \|T(x)^{-1}\|$ is locally bounded, and $\sH_1$ is separable.
Then $\sH_2$ is separable (and hence 5.3 applies).

\demo{Proof} A partition of unity argument gives a reduction to the case $\|T(x)^{-1}\|\leq 1$, $\forall x$.
Let $e_\alpha$ be as in the proof of 5.5 and let $f_k,\ k=1,2,\ldots$ be non--vanishing sections of $\sH_{1|X_k},X_k$ closed in $X$, such that $\{f_k(x)\}$ is dense in $\sH_1(x)$, $\forall x$.
For each $k$ there is an open cover $\sU_k$ of $X_k$ and an $\alpha(U)\in D$, $ \forall \ U\in \sU_k$,  such that $(e_{\alpha(U)}(x) Tf_k(x),Tf_k(x))>{1\over 2}(f_k(x),f_k(x))$, $\forall \ x\in U$.
By paracompactness, there is a refinement of $\sU_k$ of the form $\bigcup\limits^\infty_{j=1}\cV_{k,j}$, where each $\cV_{k,j}$ is a discrete (in $X$) family of sets.  (Recall that discrete means that each point of $X$ has a neighborhood that intersects at most one member of $\cV_{k,j}$.)
Then there are sections $g_{k,j}$ of $\sK(\sH_2)$ such that $0\leq g_{k,j}(x)\leq 1$ and $(g_{k,j}(x) Tf_k(x),Tf_k(x))>{1\over 2}(f_k(x),f_k(x))$, $\forall \ x\in X_k\cap (\cup\cV_{k,j})(g_{k,j}$ is made by using a different $e_\alpha$ for each $V\in\cV_{k,j}$).
Then $\forall x\in X$, $\bigvee\limits_{k,j}$ range $g_{k,j}(x)=\sH_2(x)$.
The conclusion follows with the help of 5.4(ii) (or otherwise).
\enddemo

Let $H$ be a separable infinite dimensional Hilbert space, 
$$
E=B(H)\text{ with the strong operator topology,}
$$
$$
 B=B(H)_{sa}\text{ with the weak operator topology, and}
$$ 
$$
P=B(H)_+\text{ with the weak operator topology.}
$$

There are maps $q\colon E\to P$ and $r\colon E\to B$ defined by $q(T)=T^*T$ and $r(T)={T+T^*\over 2}$.
Special cases of the positive results 4.8, 4.9 and 5.3 have to do with lifting problems for $q$, and the negative result 4.22 has to do with a lifting problem for $r$.
Therefore, it is interesting to consider $q$ and $r$ from the point of view of fibrations. In 
5.8 and 5.10 below we show that this is not totally unreasonable (it may seem that the topology of $E$ is too strong relative to that of $P$ and $B$).

{\bf 5.8.\ Proposition}.
The map $r$ is open.

\demo{Proof} Since $r$ is real--linear, it is enough to check that it is open at $0\in E$.
Let $F$ be a finite dimensional subspace and represent $h\in B$ by $\bmatrix a&b\\ b^*& c\endbmatrix$, relative to $H=F\oplus F^\perp$.
We show $\|a\|<\epsilon\Rightarrow h=r(T)$ with $\|T_{|F}\|<\epsilon$.
This is trivial, since we may take $T=h+\bmatrix\format\r&\quad\r\\ 0&b\\ -b^*&0\endbmatrix$.
\enddemo

{\bf 5.9.\ Lemma}.
Let $F$ be a finite dimensional subspace of $H$ and $0\leq h=\bmatrix h_1&r\\ r^*&h_2\endbmatrix$, relative to $H=F\oplus F^\perp$.
If $\bmatrix a\\ c\endbmatrix$ is such that $a^*a+c^*c=h_1$ and $h_1$ is invertible, then there is $T=\bmatrix a&b\\ c&d\endbmatrix$ such that $T^*T=h$.

\demo{Proof}  Choose $s$ with $\|s\|\leq 1$ such that  $r=h_1^{1/2} sh_2^{1/2}$,  and let  $\bmatrix a\\ c\endbmatrix=uh_1^{1/2}$ (polar decomposition).  (The existence of $s$ follows from a criterion for the positivity of $2\times 2$ matrices of operators which may be folklore.  One way to prove it is to first do the special case $h_1 =h_2 =1$ by using a polar decomposition of $r$, and then reduce to this case by considering

$$
\bmatrix (h_1 +\epsilon)^{-1/2})&0\\ 0&(h_2 +\epsilon)^{-1/2})\endbmatrix \bmatrix h_1 +\epsilon&r\\ r^* &h_2 +\epsilon \endbmatrix \bmatrix (h_1 +\epsilon)^{-1/2})&0\\ 0&(h_2 +\epsilon)^{-1/2})\endbmatrix .
$$

\noindent
The fact that $||(h_1+\epsilon)^{-1/2} r(h_2+\epsilon)^{-1/2})|| \le 1$, $\forall \epsilon >0$, implies the existence of $s$ by a weak compactness argument.)  
We seek $\bmatrix b\\ d\endbmatrix=vh_2^{1/2}$ with $v$ an isometry from $F^\perp$ into $H$ such that $u^*v=s$.
If we write $H=$ range $u\oplus$(range $u)^\perp$, the first component of $v$ is uniquely determined, and there is no difficulty choosing the second to make $v$ an isometry.
\enddemo

{\bf 5.10.\ Proposition}.
The map $q$ is open.

\demo{Proof} Let $F$ be a finite dimensional subspace of $H$ and $T_0=\bmatrix a_0&b_0\\ c_0&d_0\endbmatrix\in E$ (same notation as above).
We need:\ $\forall \ \epsilon > 0$, $\exists \ \delta>0$ such that if $h=\bmatrix h_1&r\\ r^*&h_2\endbmatrix\geq 0$ and $\|h_1-a_0^*a_0-c_0^*c_0\|<\delta$, then $\exists \ T=\bmatrix a&b\\ c&d\endbmatrix\in E$ with $\|\bmatrix a\\ c\endbmatrix - \bmatrix a_0\\ c_0\endbmatrix \|<\epsilon$ and $T^*T=h$.
We may assume $a_0^* a_0+c_0^* c_0$ invertible.
In fact, otherwise let $p$ be the projection on the initial space of $\bmatrix a_0\\ c_0\endbmatrix$.
If the invertible case is known, we can find $\delta$ such that $\delta^{1/2} < \epsilon/2$ and $||ph_1p - a_0^*a_0 -c_0^*c_0||< \delta$ implies the existence of $T$ with $T^*T =h$ and $||\bmatrix ap\\ cp\endbmatrix - \bmatrix a_0\\ c_0\endbmatrix|| <\epsilon/2$.  Then, given the stronger condition $||h_1 - a_0^*a_0 -c_0^*c_0||< \delta$, we have $||\bmatrix a(1-p)\\ c(1-p)\endbmatrix||<\delta^{1/2}$, whence $\|\bmatrix a\\ c\endbmatrix-\bmatrix a_0\\ c_0\endbmatrix\| < \delta^{1/2}+\|\bmatrix ap\\ cp\endbmatrix-\bmatrix a_0\\ c_0\endbmatrix \|< \epsilon$.
Now, assuming the invertibility, in view of 5.9 it is enough to find $a,c$.
Write $\bmatrix a_0\\ c_0\endbmatrix=u_0(a_0^* a_0+c_0^*c_0)^{1/2}$ (polar decomposition) and $\bmatrix a\\ c\endbmatrix=u_0 h_1^{1/2}$.
\enddemo

{\bf 5.11.\ Theorem}.
There is a function $s\colon P\to E$ such that 

(i)\ $q\circ s$ = id

(ii)\ $s|P'$ is continuous, $\forall$ bounded $P'\subset P$

(iii)\ $h$ one--one $\Rightarrow s(h)$ is one--one with dense range

(iv) $h$ invertible $\Rightarrow s(h)$ invertible.

\demo{Proof} (cf.~proof of 4.9).
Let $e$ be a strictly positive element of $\sK$.
For $h\in P$ define $s(h)$ as the unique $T\in E$ such that $Te=(ehe)^{1/2}$.
If $(h_\alpha)$ is bounded and $h_\alpha\to h$, weakly, then $eh_\alpha e\to ehe$, in norm $\Rightarrow s(h_\alpha) e\to s(h)e$ in norm $\Rightarrow s(h_\alpha)\to s(h)$, strongly, since $(s(h_\alpha))$ is bounded.
\enddemo

{\bf 5.12.\ Corollary}.
For any locally bounded continuous map $f\colon X\to P$, there is a continuous lifting $\tilde f\colon X\to E$ such that $\tilde f(x)$ is one--one with dense range or invertible whenever $f(x)$ is one--one or invertible.

Of course it is possible for $T$ to have dense range even if $q(T)=T^*T$ fails to be one--one.
The result on liftings of this type is weaker.

{\bf 5.13.\ Theorem}.
If $X$ is a finite dimensional paracompact space, if $f\colon X\to P$ is locally bounded and continuous, and if $f(x)$ has infinite rank, $\forall \ x\in X$, then there is a continuous lifting $\tilde f\colon X\to E$ such that $\tilde f(x)$ has dense range, $\forall \ x\in X$.

\demo{Proof} There is an analogue to the construction of $X_h$ in \S4.
Let $\sH_1$ be the trivial continuous field with $\sH_1(x)=H$, $\forall \ x\in X$.
For each $x\in X$, $f(x)$ gives a new semi--definite inner product on $H$.
Let $\sH_2(x)$ be the Hausdorff completion of $H$ relative to this inner product.
Then $\sH_2$ is a continuous field of Hilbert spaces, each infinite dimensional, and the ``identity map'' from $H$ to $\sH_2(x)$ gives an $L$--map $T$ such that $T(x)$ has dense range and $T(x)^*T(x)=f(x)$, $\forall \ x\in X$.
By 5.6, $\sH_2$ is separable, and hence Theorem 5, p.~260, of [21] implies that $\sH_2$ is trivial.
If $U$ is a unitary isomorphism of $\sH_1$ and $\sH_2$, then set $\tilde f(x)=U(x)^* T(x)$.
\enddemo

5.13 certainly fails if $X$ is not finite dimensional.
In fact, by [21], there is a continuous map $f\colon X\to$ the strong Grassmanian, such that $f(x)$ has infinite rank, $\forall x$, and the continuous field defined by $f$ is non--trivial, where $X$ is the product of countably many 2--spheres.  If $\tilde f$ is a continuous lifting of $f$, then from $\tilde f(x)^*\tilde f(x) =f(x)$ it follows that $\tilde f(x)$ is a partial isometry whose initial space is the range of $f(x)$.  If $\tilde f(x)$ has dense range, then it is surjective, and it sets up an isomorphism of the continuous field defined by $f$ with a trivial continuous field.

Let $E_1=\{T\in E\colon T$ is one--one with dense range$\}$.
$E_2=\{T\in E\colon T$ has dense range$\}$, and $G=\{T\in E\colon T^{-1}$ exists$\}$.
Consider $q_1=q_{|E_1}\colon E_1\to \{h\in P\colon h$ is one--one$\}$, $q_2=q_{|E_2}\colon E_2\to \{h\in P\colon h$ has infinite rank$\}$, and $q_0=q_{|G}\colon G\to \{h\in P\colon h^{-1}$ exists$\}$.
To consider whether $q, q_0,q_1$, or $q_2$ is really a fibration, we should be able to compare two liftings of the same map.
Let $\sU$ be the unitary group of $H$ with the strong topology.

{\bf 5.14.\ Lemma}.
If $X$ is any space and $\tilde f_1,\tilde f_2\colon X\to E_2$ are continuous maps such that $q\circ\tilde f_1=q\circ \tilde f_2$, then there is a continuous $g\colon X\to\sU$ such that $\tilde f_2(x)=g(x)\tilde f_1 (x)$, $\forall \ x\in X$.

\demo{Proof} $g(x)$ is unique of course since $\tilde f_i(x)$ has dense range, and it is routine to check the continuity of $g$.
\enddemo

$q_0$ and $q_1$ are still not quite fibrations because of the boundedness conditions in 5.11 and 5.12.
To remedy this, we  modify the topologies on $E$ and $P$ to the strongest topologies which agree with the original ones on bounded sets.
The reader can check that 5.8 and 5.10 remain true.
Then 5.11 and 5.14 give that the modified $q_0$ and $q_1$ are trivial fibrations with fiber $\sU$, and 5.13 and 5.14 give that $q_2$ is a Serre fibration.
Since the boundedness conditions will be satisfied in applications, we are not formalizing the modifications.
We do not know whether $q$ is a fibration.

As to $r$, if $r$ were any kind of fibration, maps of $[0,1]$ to $B$ would have to lift to $E$.
This is contradicted by 4.22, since $QM(C([0,1])\otimes\sK)$ can be identified with the set of continuous maps from $[0,1]$ to $B$ and $LM(C([0,1])\otimes\sK)$ with the set of continuous maps from $[0,1]$ to to $E$.  Since $QM(C([0,1])\otimes\sK)\ne LM(C([0,1])\otimes\sK)+RM(C([0,1])\otimes\sK)$, there is a self-adjoint $S$ in $QM(C([0,1])\otimes\sK)$ such that $S\ne T+T^*$ with $T\in LM(C([0,1])\otimes\sK)$.
On the other hand, if $X=\{0\}\cup \{{1\over n}\colon n=1,2,\ldots\}\subset [0,1]$ then maps from $X$ to $B$ do lift to $E$, by a proof similar to that of 5.8.

Since $r$ looks like the differential of $q_0$, our intuition from Lie theory seems to contradict the idea that $q_0$ should be better behaved than $r$.
$G$ is not quite a topological group, but this is not the problem.
As will be shown in 6.4, although maps into $q_0(G)$ lift to $G$, maps (from $[0,1]$, say) to $q_0(G)$ which are close to 1 in norm do not lift to maps to $G$ which are close to 1 in norm.

\subhead  \S6.\ Counter-examples\endsubhead

\example{6.1} We present the example of Choi and Christensen [12] from the point of view of 4.8 and 4.9.
Thus we show that there is a $C^*$--algebra $C$ and invertible positive $h\in QM(C)$ such that $h\neq T^*T$ for any $T\in LM(C)$.
We also show that $QM(C)\neq$ [span$(RM(C)\cdot LM(C))]^-$.

Let $\pi\colon B(H)\to B(H)/\sK=Q$ be the quotient map.
Let $A,B\subset Q$ be $C^*$--subalgebras such that $A\cdot B=0$ and there does not exist $s\in Q$ with $As=(1-s)B=0$ ([12]).
Let $C=(\bmatrix c_{11}&c_{12}\\ c_{21}& c_{22}\endbmatrix\in B(H)\otimes M_2\colon \pi(c_{11})\in A$, $\pi(c_{22})\in B, c_{12}, c_{21}\in\sK)$.
Then $QM(C)$ and $LM(C)$ can be identified with their images in $B(H)\otimes M_2$, although we are using a non--universal representation of $C$.
So $T=\bmatrix t_{11}&t_{12}\\ t_{21}&t_{22}\endbmatrix$ is a quasi--multiplier if and only if $A\pi(t_{11})A\subset A$, $B\pi (t_{22})B\subset B$, and $A\pi (t_{12})B=B\pi (t_{21})A=0$.
In particular, any scalar matrix is a quasi--multiplier.
And $T$ is a left multiplier if and only if $\pi(t_{11})A\subset A$, $\pi(t_{22})B\subset B$, and $\pi(t_{12})B=\pi(t_{21})A=0$.
Let $L_0=\{t\in B(H)\colon \pi(t)B=0\}$ and $R_0=\{t\in B(H)\colon A\pi(t)=0\}$.
Then we see that $T\in RM(C)\cdot LM(C)\Rightarrow t_{12}\in L_0+R_0$.

Now let $h=\bmatrix 1&\epsilon\\ \epsilon&1\endbmatrix$, $\epsilon>0$, which is an invertible positive quasi--multiplier.
The conditions on $A,B$ yield that $1\not\in L_0+R_0$, and hence $h\not\in \text{ span}(RM(C)\cdot (LM(C))$.
It is possible to choose $A,B$ such that $1\not\in (L_0+R_0)^-$ (for example, this is true for example (b) of [12]), and then $h\not\in [\text{span}(RM(C)\cdot LM(C))]^-$.
\endexample

\example{6.2} Let $H$ be a non--separable Hilbert space, $Y$ a suitable compact Hausdorff space (as specified below), and $B=C(Y)\otimes\sK(H)$.
We will show that $\forall \ \epsilon > 0$, $\exists \ 0\leq h\in QM(B)$ such that $\|1-h\|\leq\epsilon$ and $h\neq T^*T$, $\forall$ invertible $T\in LM(B)$.
If $X_h$ is the right Hilbert $B$--module constructed from $h$ as in \S4 and $A=\sK(X_h)$, then further $A$ is not isomorphic to $B$.
Thus we have a counterexample to (P1) as well as (P2).
The example of Choi and Christensen [12] also accomplishes the above.
Since $A$ and $B$ are derived from continuous fields of Hilbert spaces (in particular they are of continuous trace), 6.2 also gives a counter--example to question 2$^\circ$, p.~265 of [21] (cf.~proof of 5.13).

Let $P_\epsilon=\{p\in B(H)\colon 1-\epsilon\leq p\leq 1\}$ with the weak operator topology, $0<\epsilon<1$.
Let $Y= P_{1/2}$, which is homeomorphic to $P_\epsilon$ for all $\epsilon$. 
Identify $QM(B)$ with the set of weakly continuous functions from $Y$ to $B(H)$ and $LM(B)$ with the set of strongly continuous functions.
Let $h$ be a homeomorphism from $Y$ to $P_\epsilon$.

\medskip
{\bf Lemma}.\ If $T$ is invertible in $LM(B)$, then for any separable $H_0\subset H$, there is a separable subspace $H_1\supset H_0$ such that $T(y)H_1$, $T(y)^{-1} H_1\subset H_1$, $\forall \ y\in Y$.

\demo{Proof} Let $v_1,v_2,\ldots$ be a dense sequence in $H_0$.
For each $k$, $T(Y)v_k$ and $T^{-1}(Y)v_k$ are compact subsets of $H$, hence separable.
(Here we use 4.1.)
Thus there is a separable subspace $H'_1\supset H_0$ such that $T(Y)H_0$, $T^{-1}(Y) H_0\subset H'_1$.
Recursively we construct separable subspaces $H_0\subset H'_1\subset H'_2\subset\ldots$ such that $T(Y)H'_n$, $T^{-1}(Y)H'_n\subset H'_{n+1}$.
Then take $H_1=(\bigcup\limits^\infty_{n=1} H'_n)^-$.
\enddemo

Now suppose there were an invertible $T\in LM(B)$ such that $T^*T=h$.
Choose an infinite dimensional $H_1\subset H$ satisfying the lemma and represent operators by $2\times 2$ matrices relative to $H=H_1\oplus H_1^\perp$.
Let $Y_0=h^{-1} (\{\bmatrix 1-{\epsilon\over 2}&r\\ r^*&1-{\epsilon\over 2}\endbmatrix\colon \|r\|\leq {\epsilon\over 2}\})$.
Let $\Delta$ be the ball of radius ${\epsilon\over 2}$ in $B(H_1^\perp,H_1)$, with the weak operator topology.
Then since $h$ maps onto $P_\epsilon$, $\Delta$ may be regarded as a homeomorphic image of $Y_0$.
Now $T_{|Y_0}=\bmatrix a&b\\ 0&c\endbmatrix$, relative to $H=H_1 \oplus H_1^\perp$, where $a^*a=1-{\epsilon\over 2}$, $b^*b+c^*c=1-{\epsilon\over 2}$, and $a^*b=r$.
Since $H_1$ is invariant under $T^{-1}$, $a$ is invertible and $(1-{\epsilon\over 2})^{-1/2}a$ is unitary.
Since the map $u\mapsto u^*$ is strongly continuous for unitary $u$, we see that $a^*$ is a strongly continuous function on $Y_0$.
This means that $r=a^*b$ is strongly continuous; hence the strong topology is the same as the weak topology on $\Delta$, a contradiction.

To see that $A$ is not isomorphic to $B$, note that $A$ comes from a continuous field of Hilbert spaces $\sH_2$ on $Y$ (as in the proof of 5.13).
We have proved that $\sH_2$ is not the trivial field $\sH_1$.
Since there are automorphisms of $B$ inducing any self--homeomorphism of $Y$, if $A$ and $B$ were isomorphic, there would be an isomorphism $\varphi\colon A\to B$ over $C(Y)$.
By Theorem 9, p.~272 of [21] $\varphi$ is induced by an isomorphism of $\sH_2$ with $\sH_1\otimes\ell$, for some line bundle $\ell$ on $Y$.
But $\sH_1\otimes\ell$ is the direct sum of $\ell$ with itself uncountably many times, and Corollary 3, p.~260 of [21] implies $\sH_1\otimes\ell$ is trivial.
Thus $\sH_2$ is trivial, a contradiction.
\endexample

\example{6.3} It is not difficult to find an example of two continuous fields of Hilbert spaces, $\sH_1$, and $\sH_2$, on a compact Hausdorff space $X$, such that $\sH_1$ is separable, $\sH_2$ is not separable, and there exists an $R$--map (Def 5.1) $T$ from $\sH_1$ to $\sH_2$ such that $T(x)$ is one--one with dense range, $\forall \ x\in X$.
Let $X$ be any space with a point $p$ such that $\{p\}$ is not a $G_\delta$.
Let $H$ be a separable infinite dimensional Hilbert space and $\sH_1$ the trivial field with fiber $H$.
Let $\ell'$ be the subfield of the trivial line bundle such that $\ell'(x)=\bC$ for $x\neq p$ and $\ell'(p)=0$.
Thus the continuous sections of $\ell'$ are the continuous scalar--valued functions on $X$ vanishing at $p$.
Let $\sH_2=\sH_1\oplus\ell'$.
By [21, Prop.~13 p.~242] $\sH_2$ is not separable.
Let $S\in B(H)$ be one--one with dense but not closed range.
Let $v\in H\setminus$ range $S^*$.
Define $T(x)\colon\sH_1(x)\to\sH_2(x)$ by
$$
T(x)u=\cases Su,&x=p\\ Su\oplus(u,v),&x\neq p.\endcases
$$
\endexample

{\bf 6.4.\ Theorem}.
Let $A$ be a $C^*$--algebra such that $QM(A)\neq LM(A)+RM(A)$.
Then $\exists \epsilon>0$ such that $\forall \ \delta>0$, $\exists h\in QM(A)$ such that $1-\delta\leq h\leq 1+\delta$ and $\not\exists \ T\in LM(A)$ with $\|T-1\|<\epsilon$ and $T^*T=h$.

\demo{Proof} First note that if $\|T-1\|$ is sufficiently small then $T=e^s$ where $\|s\|$ is small enough for convergence of the Campbell--Baker--Hausdorff formula.
Also $\|s\|\sim \|T-1\|$ and $s$ is again in $LM(A)$, since $LM(A)$ is a Banach algebra.

{\bf Lemma}.
Write $h=T^*T$ with $T=e^s$, as above.
Then $h-1=s+s^*+O(\|h-1\|\cdot\|T-1\|)$ as $T\to 1$ in norm.

\demo{Proof} Write $s=k+ij$, $k,j$ self--adjoint and $h=e^{h_0}$ $(h_0=h-1-{(h-1)^2\over 2}+\ldots)$.
It is not claimed that $h_0\in QM(A)$.
Then $e^{h_0}=e^{s^*} e^s\Rightarrow h_0 = 2k+a_1[k,j]+\ldots$ (Campbell--Baker--Hausdorff plus $[s^*,s]=2i[k,j]$), where all of the omitted terms involve $[k,j]$.
Thus $h_0=2k+O(\|k\|\cdot\|j\|)=2k+O(\|k\|\cdot\|s\|)$.
Now $h=1+h_0+O(\|h_0\|^2)$ and $\|h-1\|\sim \|h_0\| \sim 2\|k\|$.
Since $\|s\|\geq \|k\|$, $\|h_0\|^2=O(\|k\|\cdot\|s\|)$, and $h-1=2k+O(\|k\|\cdot\|s\|)$ or $h-1=s+s^*+O(\|h-1\|\cdot\|T-1\|)$.
\enddemo

Thus $\exists \ M$, $\epsilon>0$ such that $\|h-1-s-s^*\|\leq M\|T-1\|\cdot\|h-1\|$ and $\|s\|\leq 2\|T-1\|$ whenever $\|T-1\|\leq\epsilon$.
We may assume $\epsilon M\leq {1\over 2}$.
Now assume that there is a $\delta\in (0,1)$ such that $h\in QM(A)$ and $1-\delta\leq h\leq 1+\delta\Rightarrow h=T^*T$ with $T\in LM(A)$ and $\|T-1\|<\epsilon$.
Let $h_1\in QM(A)_{sa}$ and take $h=1+{\delta\over \|h_1\|}\ h_1$.
Then $\exists \ s\in LM(A)$ such that $\|s\|<2\epsilon$ and ${\delta\over \|h_1\|}\ h_1=2\text{Re}\ s+h'_2$ with $\|h'_2\|\leq {1\over 2}\delta$.
Thus $h_1=$ \text{Re} $s_1+h_2$, where $s_1\in LM(A)$, $\|s_1\|\leq {4\epsilon\over\delta} \|h_1\|$, and $\|h_2\|\leq{1\over 2}\|h_1\|$.
By a familiar recursive procedure, we find $h_1= \text{Re } R$, $R\in LM(A)$ (and $\|R\|\leq {8\epsilon\over\delta} \|h_1\|$).
This means $QM(A)=LM(A)+RM(A)$, a contradiction.
\enddemo

\example{Remarks}

(i)\ The $\epsilon$ in the theorem is a universal constant.
To estimate it, it would be necessary only to be more careful in proving the lemma.

(ii)\ The converse of the theorem is also true.
Let $G_n=\{\text{invertible }T\in LM(A)\}$ and $P_n=\{\text{invertible }h\in QM(A)_+\}$ both with the norm topology, and $q_n\colon G_n\to P_n$ by $q_n(T)=T^*T$.
Note that $G_n$ is a topological group, $G_n$ acts on $P_n$ from the right by $(h,T)\mapsto T^*hT$, and $q_n$ is just the map from $G_n$ to the orbit of $1\in P_n$.
If $QM(A)=LM(A)+RM(A)$, then $q_n$ is open and surjective, even if $A$ is not $\sigma$--unital.

We wish to derive from the theorem an example of the following:\ A ``nice'' $C^*$--algebra $A$ and a complete order automorphism $\varphi$ of $A$ such that $\varphi$ is almost isometric $(\|\varphi\|-1$ and $\|\varphi^{-1}\|-1$ are small) but there is no automorphism $\theta$ of $A$ such that $\|\theta-\varphi\|$ is small.

The basic method of attempting this is clear.
Let $A$ be a $\sigma$--unital $C^*$--algebra such that $QM(A)\neq LM(A)+RM(A)$ (see 4.22) and $0\leq h\in QM(A)$ such that $\|h-1\|$ is small.
By 4.8, there is an invertible $T\in LM(A)$ such that $T^*T=h$.
If $\varphi(a)=TaT^*$, then $\varphi$ is almost isometric.
If $\theta$ is an automorphism such that $\|\theta^{-1}-\varphi\|$ is small, then $\theta\circ\varphi$ is a complete order automorphism close to id.
If we can prove that this implies $h=T_0^*T_0$ with $T_0\in LM(A)$ and $\|T_0-1\|$ small, we will contradict the theorem and establish (non--constructively) the desired example.
The most obvious way to do this is first to prove that $\theta$ is inner.
Then by replacing $T$ with $UT$ for some unitary $U\in M(A)$, we may assume $\theta=id$.
The problem then becomes to deduce from $\|\varphi-id\|$ small that there is a unitary $V\in M(A)$ with $\|T-V^*\|$ small.
We sketch below how to carry this out if $A$ is either simple or continuous trace.
The success in these cases seems to involve the fact that bounded derivations are inner for these algebras, but it is possible that better methods applicable to more general algebras exist.
(It is not really necessary to prove $\theta$ inner.
It would be sufficient to prove $\theta$ is close to an inner automorphism.)

If $A$ is simple, by Kishimoto [30] $\theta$ will be inner if there is $S$ and $\delta>0$ such that $\|\theta(a^*)Sa\|\geq \delta\|a\|^2$, $\forall \ a\in A$.
In [30], $S$ was required to be in $M(A)$, but in fact it is sufficient that $S\in QM(A)$.
(The author has been intending to write a paper generalizing Kishimoto's results for several years.
That $S\in QM(A)$ is sufficient is minor.)
If $\|\theta\circ\varphi-id\|$ is small, it is easy to see that $S=\theta(T^*)$ will suffice in Kishimoto's theorem.
Now assume $\theta=1$.
Write $T=Uh^{1/2}$ (polar decomposition).
If $\pi$ is an irreducible representation of $A$ (extended to $A^{**}$), then $\|\pi(T)\cdot\pi(T^*)-id\|$ is small on $\pi(A)''=B(H)$, in view of the Kaplansky density theorem.
Since $\|\pi(h^{1/2})-1\|$ is small, [34, 8.7.5] $\Rightarrow \exists \ \lambda$ such that $|\lambda|=1$ and $\|\pi(U)-\lambda\|$ is small, and hence $\|\pi(T)-\lambda\|$ is small.
Since any faithful representation of $A$ is isometric on $LM(A)$ (even $QM(A))$, $\|T-\lambda\|$ is small and we may take $T_0=\overline\lambda T$.

If $A$ is continuous trace, $A$ comes from a continuous field of elementary $C^*$--algebras on $X=\hat A$.
Since $A$ is $\sigma$--unital, $X$ is $\sigma$--compact and hence paracompact.
If $\theta$ acts non--trivially on $X$, then $\exists \ 0\neq a\in A$ such that $\theta(a)$ and $a$ have disjoint supports; and this contradicts $\|\theta\cdot\varphi-id\|$ small.
Then using Theorem 9, p.~272, of [21], we can see that $\theta$ is locally inner.
Say $\theta$ is Ad $U$, where $U$ is a locally defined unitary function.
Then, as above, for each $x\in X$ there is $\lambda(x)$ such that $|\lambda(x)|=1$ and $\|\overline{\lambda(x)} U(x) T(x)-1\|$ is small.
Since $\lambda(x)$ is not uniquely defined, we wish to normalize it so that it becomes unique and continuous.
To do this, we note that there is a continuous global function $\rho$ such that each $\rho(x)$ is a state (not pure) on the elementary $C^*$--algebra $A(x)$.
In fact, locally $A$ comes from a continuous field of Hilbert spaces possessing a unit vector section.
Thus, locally we have continuous fields of pure states, and we can patch these with a partition of unity.
From the construction of $\rho$, $\rho(x)(U(x)T(x))$ is continuous; and we normalize $\lambda$ by the requirement $\rho(x)(\overline{\lambda(x)}U(x)T(x))>0$.
Then if $V(x)=\overline\lambda(x) U(x)$, $V$ is actually globally defined and gives a unitary in $M(A)$.
Take $T_0=VT$.
\endexample

\example{6.5} Let $A$ be a $C^*$--algebra such that $QM(A)\neq LM(A)+RM(A)$, and let $B=A\otimes M_2$.
We show that there is an invertible $T\in QM(B)$ such that $T^{-1}\in QM(B)$ also and $T\neq RL$ for any invertible $R\in RM(B)$, $L\in LM(B)$.
Let $S\in QM(A)\backslash (LM(A)+RM(A))$ and take $T=\bmatrix 1&S\\ 0&1\endbmatrix$.
If $T=RL$, as above, write $R=\bmatrix a&b\\ c&d\endbmatrix$, $L^{-1}=\bmatrix e&f\\ g&h\endbmatrix$, so that $\bmatrix a&b\\ c&d\endbmatrix=\bmatrix 1&S\\ 0&1\endbmatrix \bmatrix e&f\\ g&h\endbmatrix$, $a,b,c,d\in RM(A)$, $e,f,g,h\in LM(A)$.
Compute $a=e+Sg$, $b=f+Sh$, $c=g$, and $d=h$.
Thus, $g,h\in LM(A)\cap RM(A)=M(A)$.
Let $H$ be the Hilbert space of the universal representation of $A$, and consider $(g\ h)$ as an operator from $H\oplus H$ to $H$.
Then $(g\ h)$ is surjective, since it is a row of an invertible matrix.
This implies $gg^*+hh^*$ is invertible in $B(H)$.
But since $gg^*+hh^*\in M(A)\subset B(H)$, $(gg^*+hh^*)^{-1}\in M(A)$.
Hence $\exists\ r,s\in M(A)$ such that $gr+hs=1$.
Now $Sg=a-e\in LM(A)+RM(A)$ and $Sh=b-f\in LM(A)+RM(A)$.
Thus $S=(Sg)r+(Sh)s\in [LM(A)+RM(A)]\cdot M(A)\subset LM(A)+RM(A)$, a contradiction.
\endexample

\example{6.6} We give examples of nice $C^*$--algebras for which $\exists T\in QM(A)$ not of the form $R\cdot L$, $R\in RM(A)$, $L\in LM(A)$.

(a)\ Let 
$$
A=\left\{ T=\bmatrix t_{11}&t_{12}\\ t_{21}&t_{22}\endbmatrix\in B(H)\otimes M_2\colon t_{11}-\lambda , t_{12}, t_{21}, t_{22}\in\sK, \text {for some } \lambda \in\bC\right\}
$$
Then $QM(A)=(T\colon t_{11}\in\bC+\sK\}$ and $LM(A)=\{T\colon t_{11}\in\bC+\sK, t_{21}\in\sK\}$.
Then $\bmatrix 0&1\\ 1&0\endbmatrix$ is an element of $QM(A)$ not of the form $R\cdot L$.

(b)\ Let $A$ be the algebra of convergent sequences in $M_2$ with limit of the form $\bmatrix *&0\\ 0&0\endbmatrix$.
Then 
$$
QM(A)=\left\{ \matrix (x_n)^\infty_{n=1}\colon x_\infty\in\bC,\ x_n\in M_2,\ \|x_n\|\text{ bounded},\\
(x_n)_{11}\to x_\infty\endmatrix \right\}
$$
and
$$
LM(A)=\left\{ \matrix (x_n)\colon x_\infty\in\bC,\ x_n\in M_2,\ \|x_n\|\text{ bounded},\\
(x_n)_{11}\to x_\infty, (x_n)_{21}\to 0\endmatrix \right\}.
$$
Let $T\in QM(A)$ be given by $(x_n)$ where $x_\infty=0$ and $x_n=\bmatrix 0&1\\ 1&0\endbmatrix$, $n=1,2,\ldots$
Then $T$ is not of the form $R\cdot L$.
\endexample

\example{6.7} It might be thought that the hypotheses of 3.7 and 5.7 could be weakened by requiring that $T$ be surjective instead of invertible.
This can easily be seen to be wrong, even if $T\in RM(X)$ (if $T\in LM(X)$ 3.6 would apply).
Let $B$ be a $\sigma$--unital $C^*$--algebra with a non--$\sigma$--unital hereditary $C^*$--subalgebra $A$, such that $A$ generates $B$ as an ideal.
Let $X=(AB)^-$.
If $p$ is the open projection for $A$, then take $T=p$ regarded as an element of $X^{**}$.

To see that this situation is widespread, let $A_0$ be any non--$\sigma$--unital $C^*$--algebra and $\tilde A_0$ the result of adjoining a unit.
Take $B=\tilde A_0\otimes M_2$ and $A=\bmatrix \tilde A_0&A_0\\ A_0&A_0\endbmatrix\subset B$.
To see that this can occur with continuous fields of Hilbert spaces, one could use the same $\sH_1$ and $\sH_2$ that were used in 6.3.
\endexample

\subhead \S7.\ Remarks and questions\endsubhead

\example{7.1} The positive results on Takesaki's question and example 6.4 can be viewed from the point of view of perturbations of $C^*$--algebras (see [13], [14], [25], [26], [29], [36], [37] for example).
Perturbations of $C^*$--algebras have been considered in several different ways (not all of which are mentioned below).
One is to ask whether a linear almost isometric map (which could be assumed positive or completely positive) between $C^*$--algebras is close to an isometry (which could be required to be positive or completely positive, in which case it would have to be a Jordan isomorphism or isomorphism).
Another is to consider $C^*$--subalgebras $A$ and $B$ of $B(H)$ which are close in the sense that the Hausdorff distance between their unit balls is small, and ask whether there is an isomorphism $\theta\colon A\to B$ such that $\|\theta-id_{|A}\|$ is small.
In both cases one could weaken the question and simply ask whether the two algebras have to be isometric or isomorphic.

We now explain the intersection of our results with the second problem mentioned above (the relation with the first is clear).
If $A$ and $B$ are hereditary $C^*$--subalgebras of $C$ with open projections $p$ and $q$, and if $p\in M(C)$, then $pbp=pqbqp\in A$, $\forall \ b\in B$.
Thus if $\|p-q\|$ is small, $A$ and $B$ are close in the above sense.
(Of course we can consider that $C\subset B(H)$.)
Conversely, if $A$ and $B$ are close in the above sense, we can see that $\|p-q\|$ is small (even without assuming $p\in M(C)$):  If $(e_i)$ is an approximate identity of $B$ and if the Hausdorff distance between the two unit balls is less than $\epsilon$, then there is $a_i$ in $A$ such that $||a_i-e_i||<\epsilon$.  Since $||(1-p)a_i||=0$, then $||(1-p)e_i||<\epsilon$.  Taking weak limits, we find that $||(1-p)q||\le\epsilon$.   Similarly $||(1-q)p||\le\epsilon$.
As mentioned in Proposition 2.2, when $p\in M(C)$, $pbp=TbT^*=\varphi(b)$.
Thus the question whether $\varphi$ is close to an isomorphism is the same as the question whether there is an isomorphism such that $\|\theta-id_{|A}\|$ is small in this case.
Finally, the construction of N.T.~Shen [39] and the analysis of complete order isomorphisms given in \S4 show that every instance of the first perturbation problem, with the completely positive interpretation, does arise in this context (and is also an instance of the second perturbation problem).
For Shen's construction it is necessary that $\|T\|\leq 1$, but this can easily and harmlessly be achieved by multiplying $T$ by a number slightly less than 1.

In [26] B.E.~Johnson gave a counter--example for the second perturbation problem mentioned above with $A$ and $B$ both isomorphic to $C([0,1])\otimes\sK$.
Our 6.4 uses ideas very similar to Johnson's and gives a similar counter--example for the first perturbation problem.
So far as we know, 6.4 is not the same example as Johnson's example, and [26] does not imply a counter--example to the first perturbation problem; but we are not at all sure of this.

Finally we should mention that the examples of Choi and Christensen [12] and our 6.2 are non--separable counter--examples to both problems.
\endexample

\example{7.2. Questions, mainly for the non--$\sigma$--unital case}

(i)\ If $A$ and $B$ are hereditary $C^*$--subalgebras of $C$ with open projections $p$ and $q$ such that $\|p-q\|<1$, are $A$ and $B$ completely order isomorphic?

(ii)\ Say that $A\sim B$ if $A$ and $B$ can be embedded in $C$ as in (i), or equivalently if there is an $A-B$ imprimitivity bimodule $X$ and an invertible $T\in QM(X)$.
Is $\sim$ an equivalence relation?
If not, is there a good description of the equivalence relation it generates?

(ii)$'$\ Same question for $\overset '\to \sim$, where $A\overset '\to\sim B$ means there is a $A-B$ imprimitivity bimodule $X$ and $T\in QM(X)$ such that $T$ is one--one with dense range.

(iii)\ Is $LM(A)+RM(A)$ strictly dense in $QM(A)$ (cf.~4.16)?
Could it even be norm dense?

(iv)\ Consider algebras of the form $C(X)\otimes\sK (H)$, where $H$ is a non--separable Hilbert space.
In 6.2 we showed that (P2) fails if $X$ is ``large''.
We do not know whether 4.9 holds for those algebras.
((P2) and 4.9 would hold if $X$ has a countable dense set.)
\endexample

\example{7.3. Remark} The relation between $LM(A)$ and $QM(A)$ has something to do with triangularity.
This is seen from 4.19, 4.20, other parts of this paper, and possibly from the reader's favorite examples.
An earlier proof of (P2), for algebras with a countable approximate identity consisting of projections, made explicit use of triangularity.

\medskip
{\bf Lemma}.\ Let $h$ be an $n\times n$ matrix with operator entries such that $0< \epsilon\leq h\leq M$ for some $\epsilon, M\in\bR$.
Then there is a unique upper triangular $T$, with 1's on the main diagonal, such that $\Delta=T^* hT$ is diagonal.
Moreover $\epsilon\leq \Delta\leq M$.

The lemma is a version of the Gram--Schmidt process, and we claim no originality for it.
It is sufficient to prove the second sentence in the case $n=2$.
It is then possible to extend the lemma to infinite matrices, and the second sentence gives control on the norms $(\|T\|, \|T^{-1}\|\leq ({M\over\epsilon})^{1/2})$.
\endexample

\example{7.4.  Historical remarks} This paper is a slightly revised version of MSRI preprint no.~11211-85, from 1985.  We will now briefly discuss some related work done afterwards.

The previously cited paper [9] contains a new proof of the main result of Shen's thesis as well as some additional results.  As suggested above, Theorems 4.13 and 4.14 were inspired by [35, Theorem 10], the non-commutative Tietze extension theorem.  We later proved other results inspired by this theorem:  [5, Subsection 3.A], [5, Theorem 3.43], and [7, Theorem 3.2].  Theorem 4.15 was strengthened in [6]:  The space $\tilde\sB_0 +i\tilde\sB_0$ is actually a $C^*$-algebra.  H. Lin [31] strengthened our results on the question $QM(A)=LM(A)+RM(A)$.  In particular, [31, Theorem 6.3] solves the problem of when this is so if $A$ is separable and stable.  Christensen, Sinclair, Smith, White, and Winter [15] proves among other things that separable nuclear $C^*$-subalgebras of $B(H)$ whose unit balls are close in the Hausdorff metric are unitarily conjugate.  This version of the perturbation problem does not demand that the isomorphism be close to the identity and goes back to a problem posed in [29].
\endexample

\bye